\documentclass[11pt,letterpaper,reqno]{amsart}
\usepackage[T1]{fontenc}
\usepackage{lmodern}
\usepackage{amsmath,amssymb,amsthm,mathtools}
\usepackage{microtype}
\usepackage{xcolor} 

\def\Xint#1{\mathchoice
{\XXint\displaystyle\textstyle{#1}}%
{\XXint\textstyle\scriptstyle{#1}}%
{\XXint\scriptstyle\scriptscriptstyle{#1}}%
{\XXint\scriptscriptstyle\scriptscriptstyle{#1}}%
\!\int}
\def\XXint#1#2#3{{\setbox0=\hbox{$#1{#2#3}{\int}$}
\vcenter{\hbox{$#2#3$}}\kern-.5\wd0}}
\def\fint{\Xint-}
\usepackage[colorlinks=true,linkcolor=blue,citecolor=blue,urlcolor=blue]{hyperref}

\renewcommand{\le}{\leqslant}
\renewcommand{\ge}{\geqslant}
\renewcommand{\phi}{\varphi}
\renewcommand{\epsilon}{\varepsilon}

\newcommand{\R}{\mathbb R}
\DeclareMathOperator*{\esssup}{ess\,sup}
\DeclareMathOperator{\diver}{div}

\newtheorem{theorem}{Theorem}[section]
\newtheorem{lemma}[theorem]{Lemma}
\newtheorem{definition}[theorem]{Definition}
\newtheorem{proposition}[theorem]{Proposition}
\newtheorem{corollary}[theorem]{Corollary}
\theoremstyle{remark}

\numberwithin{equation}{section}
\allowdisplaybreaks[2]

\title[elliptic systems]{Regularity of elliptic systems with Uhlenbeck structure}
\author{Jihoon Ok}
\address{Department of Mathematics,
Institute for Mathematical and Data Science, Sogang University,
35 Baekbeom-ro, Mapo-gu, Seoul 04107, Republic of Korea}
\email{jihoonok@sogang.ac.kr}
\date{}



\begin{document}

\begin{abstract}
We study autonomous elliptic systems with Uhlenbeck structure of the form
\[
-\operatorname{div}\big(\mathfrak a(|Du|)Du\big)=0 \quad \text{in }\ \Omega\subset \mathbb R^n
\]
Under two-sided quantitative $(p,q)$-growth conditions on $\mathfrak a$, we prove gradient H\"older continuity with constants and exponent depending only on the dimensions, $p$ and $q$.
This improves the classical theory initiated by Uhlenbeck by removing the continuity requirement on $\mathfrak a'$.
The proof uses nonlinear comparison systems whose regularity follows from a scalar uniformly elliptic equation.
Under weaker $(p,q)$-growth conditions on the convex function
$\varphi(t):=\int_0^t s\mathfrak a(s)\,ds$, we obtain local $W^{1,s}$-regularity for every $s>1$,
and hence H\"older continuity of solutions with every exponent less than one. 
An additional lower growth condition on $\mathfrak a$ yields local Lipschitz continuity, even when $\mathfrak a$ is discontinuous.
\end{abstract}

\maketitle

\section{Introduction and main results}\label{sec:intro}
We study the local regularity of weak solutions to autonomous nonlinear elliptic systems with Uhlenbeck structure.
Our aim is to determine which regularity conclusions follow from growth and convexity assumptions, without separate continuity assumptions on the nonlinearity or its derivative.
We obtain gradient H\"older continuity, local $W^{1,s}$-regularity for every  $s>1$, and local Lipschitz continuity under different $(p,q)$-growth assumptions. 
Here $p$ and $q$ describe lower and upper growth rates, following the terminology of Marcellini \cite{Marcellini1989,Marcellini1991,Marcellini1993}.

In the scalar equation case, Ladyzhenskaya and Ural'tseva \cite{LU1968} established $C^{0,\gamma}$- and $C^{1,\gamma}$-regularity for quasilinear equations of $p$-Laplace type.
Fundamental contributions to gradient H\"older regularity for the $p$-Laplace equation and related variational problems were made in \cite{Evans1982,Lewis1983,DiBenedetto1983,Manfredi1988}.
Lieberman \cite{Lieberman1991} treated more general growth conditions; see also \cite{Antonini2026}.
For general autonomous systems, minimizers may fail to be $C^1$, Lipschitz, or even locally bounded \cite{Necas1977,SverakYan2002}.
See also the classical counterexamples in \cite{DeGiorgi1968,GiustiMiranda1968}.
Thus, everywhere regularity for elliptic systems requires additional structural assumptions. 

A fundamental example is the isotropic \emph{Uhlenbeck} structure, where the nonlinearity is a scalar function of $|Du|$ multiplied by $Du$.
In her landmark work \cite{Uhlenbeck1977}, Uhlenbeck established everywhere $C^{1,\gamma}$-regularity for a class of $p$-Laplacian type systems with $p\ge2$.
Subsequently, Tolksdorf \cite{Tolksdorf1983} extended this theory to a broader class of systems for all $p>1$. 
See also \cite{AcerbiFusco1989} for the subquadratic case.
Local Lipschitz estimates under nonstandard growth conditions were established in \cite{Lieberman1993,CupiniGuidorziMascolo2003,BenedettiMascolo2004,MarcelliniPapi2006}.
For systems with Uhlenbeck structure, Diening, Stroffolini
and Verde \cite{DSV2009} proved $C^{1,\gamma}$-regularity under $(p,q)$-growth, while gradient bounds and potential estimates for nonhomogeneous systems were obtained in \cite{DuzaarMingione2010,KuusiMingione2018}.
For an overview of the regularity theory for elliptic equations and systems and these counterexamples, we refer to Mingione \cite{Mingione2006}.

More precisely, after a change of notation, the systems considered by Uhlenbeck \cite{Uhlenbeck1977} take the form
\begin{equation}\label{eq:Uhlenbeck}
 -\diver\big(\mathfrak a(|Du|)Du\big)=0
 \qquad\text{in }\Omega,
\end{equation}
where $\Omega\subset\R^n$ ($n\ge2$) is open, $\mathfrak a\in C([0,\infty))\cap C^1((0,\infty))$, and $u:\Omega\to\R^N$ with $N\ge1$.
She proved that if the nonnegative function $\mathfrak a$ satisfies, for some $p\ge2$, $\mu\ge0$, $L\ge1$ and $\gamma_0>0$,
\[
 L^{-1}\big(\mu^2+t^2\big)^{(p-2)/2}
 \le \mathfrak a(t)+t\mathfrak a'(t)
 \le L\big(\mu^2+t^2\big)^{(p-2)/2}
 \qquad\text{for }t>0,
\]
and
\[
 |t_1\mathfrak a'(t_1)-t_2\mathfrak a'(t_2)|
 \le L\big(\mu^2+t_1^2+t_2^2\big)^{(p-2)/2-\gamma_0}
          |t_1^2-t_2^2|^{\gamma_0}
 \qquad\text{for }t_1,t_2>0,
\]
then every weak solution
$u\in W^{1,p}_{\mathrm{loc}}(\Omega;\R^N)$ belongs to
$C^{1,\gamma}_{\mathrm{loc}}(\Omega;\R^N)$ for some $\gamma\in(0,1)$.
This includes the $p$-Laplace system for $p\ge2$, with
$\mathfrak a(t)=t^{p-2}$ and $\mu=0$.
When $\mu=0$ and $p>2$, the above growth condition allows degeneracy at $Du=0$.
The continuity condition on $t\mathfrak a'(t)$ is essential to the classical perturbation argument, where it controls the nonlinear Taylor remainder;
see \cite[Lemma~4.4]{Uhlenbeck1977}.

For scalar equations, gradient H\"older regularity is available
under quantitative convexity and growth assumptions without
such a continuity condition.
For systems, however, this condition has played an essential
role in the classical proofs.
Our first main result removes this requirement for systems with Uhlenbeck structure.
The key step is to retain a nonlinear comparison system and prove its regularity directly, rather than compare the original nonlinearity with a linear approximation.
For the statement below, set $\varphi(t):=\int_0^t s\mathfrak a(s)\,ds$, which is an N-function (see Definition~\ref{def:N}), and define
\begin{equation}\label{eq:AV}
 V(P):=\sqrt{\frac{\varphi'(|P|)}{|P|}}\,P
 \qquad\text{and}\qquad
 \mathcal E(u,B):=\fint_B\big|V(Du)-(V(Du))_B\big|^2\,dx,
\end{equation}
where $(V(Du))_B:=\fint_B V(Du)\,dx$.

\begin{theorem}\label{thm:gradient}
Suppose $1<p\le q<\infty$ and  $\mathfrak a:(0,\infty)\to(0,\infty)$ satisfies 
\begin{equation}\label{eq:a-indices}
 \mathfrak a \in C((0,\infty))
 \quad\text{and}\quad
 \lambda^{p-2} \mathfrak a (t)\le \mathfrak a (\lambda t) \le\lambda^{q-2}\mathfrak a (t) 
 \quad\text{for all }\ \lambda\ge1\ \text{ and }\ t>0.
 \end{equation}
Let $u\in W^{1,\varphi}_{\mathrm{loc}}(\Omega;\R^N)$ be a weak solution to \eqref{eq:Uhlenbeck}.
Then there exist constants $\sigma\in(0,1)$ and $c>0$,
depending only on $n,N,p,q$, such that
\begin{equation}\label{eq:main-excess}
 \mathcal E(u,B_\rho(x_0))
 \le c\left(\frac{\rho}{r}\right)^\sigma\mathcal E(u,B_r(x_0))
\end{equation}
for every ball $B_r(x_0)\Subset\Omega$ and every $\rho \in (0,r)$. 
In particular,
\[
 u\in C^{1,\gamma}_{\mathrm{loc}}(\Omega;\R^N),
\]
where $\gamma=\sigma/\max\{2,q\}>0$, and there exists a constant $c>0$ depending only on $n,N,p,q$ such that for every ball $B_{2r}(x_0)\Subset\Omega$,
\begin{equation}\label{eq:strong-main-holder}
 \sup_{B_r(x_0)}|Du| +r^\gamma[Du]_{C^{0,\gamma}(B_r(x_0))}
 \le c\,\varphi^{-1}\left(\fint_{B_{2r}(x_0)}\varphi(|Du|)\,dx\right).
\end{equation}
\end{theorem}

To compare the assumptions and state our remaining results, we write systems with Uhlenbeck structure as
\begin{equation}\label{main_system}
 -\diver\left(\frac{\varphi'(|Du|)}{|Du|}Du\right)=0
 \quad\text{in }\ \Omega,
\end{equation}
where $\varphi:[0,\infty)\to[0,\infty)$ is an N-function and $\varphi'$ denotes
its right derivative. 
The relation with \eqref{eq:Uhlenbeck} is $\varphi'(t)=t\mathfrak a(t)$.
When $\varphi$ is classically differentiable, \eqref{main_system} is the Euler--Lagrange system of the isotropic functional
\begin{equation}\label{energy}
 v\longmapsto\int_\Omega\varphi(|Dv|)\,dx.
\end{equation}
A map $u\in W^{1,\varphi}_{\mathrm{loc}}(\Omega;\R^N)$ is a weak solution to \eqref{main_system} if
\[
 \int_\Omega\frac{\varphi'(|Du|)}{|Du|}Du:D\zeta\,dx=0
 \qquad\text{for all }\zeta\in C_c^\infty(\Omega;\R^N).
\]

The closest result to Theorem~\ref{thm:gradient} is that of Diening, Stroffolini and Verde \cite{DSV2009}, who proved gradient H\"older continuity under two structural conditions on $\varphi$.
First, $\varphi\in C^1([0,\infty))\cap C^2((0,\infty))$ satisfies
\begin{equation}\label{eq:indices}
 0<p-1\le\frac{t\varphi''(t)}{\varphi'(t)}
 \le q-1<\infty
 \quad\text{for all }t>0.
\end{equation}
Second, $\varphi''$ satisfies
\begin{equation}\label{eq:omitted}
 |\varphi''(t+s)-\varphi''(t)|
 \le c\,\varphi''(t)\,\omega\left(\frac{|s|}{t}\right)
 \quad\text{for all }t>0\text{ and }|s|<\frac t2,
\end{equation}
where $\omega(r)=r^{\gamma_1}$ for some $\gamma_1>0$.
In the linear comparison argument of \cite{DSV2009}, condition \eqref{eq:omitted} controls the difference between the derivative of the nonlinearity and its value at a fixed nonzero gradient matrix.

Theorem~\ref{thm:gradient} removes \eqref{eq:omitted} and does not require classical second differentiability of $\varphi$.
Indeed, under the relation $\varphi'(t)=t\mathfrak a(t)$, condition \eqref{eq:a-indices} is equivalent to
\begin{equation}\label{eq:flux-indices}
 \varphi'\in C([0,\infty)),
 \quad
 \lambda^{p-1}\varphi'(t)\le \varphi'(\lambda t)\le\lambda^{q-1}\varphi'(t)
 \quad\text{for all }\lambda\ge1\text{ and }t>0.
\end{equation}
or, equivalently,
\begin{equation}\label{eq:flux-indices0}
 \varphi'\in C([0,\infty)),
 \quad
 \frac{\varphi'(t)}{t^{p-1}}\ \ \text{is increasing} \ \ \text{and}\ \ \frac{\varphi'(t)}{t^{q-1}}\ \ \text{is decreasing} \ \ \text{for }\ t>0,
\end{equation}
or
\[
\varphi'\in AC_{\mathrm{loc}}((0,\infty))
\quad \text{and}\quad
p-1\le\frac{t\varphi''(t)}{\varphi'(t)}
\le q-1
\quad\text{a.e. }\ t>0.
\]
These equivalent conditions imply that $\varphi'$ is locally Lipschitz on $(0,\infty)$. If $\varphi\in C^2((0,\infty))$, they are also equivalent to \eqref{eq:indices}.

We next consider a weaker $(p,q)$-growth condition on $\varphi$, with $1<p\le q<\infty$:
\begin{equation}\label{eq:potential-scaling}
 \lambda^p\varphi(t)\le\varphi(\lambda t)\le\lambda^q\varphi(t)
 \quad\text{for all }\ \lambda\ge1\ \text{ and }\ t\ge0.
\end{equation}
Since $\varphi$ is an N-function, this is equivalent to
\begin{equation}\label{eq:weak-indices}
 1<p\le\frac{t\varphi'(t)}{\varphi(t)}
 \le q<\infty
 \quad\text{for all }t>0.
\end{equation}
Indeed, taking the right derivative in $\lambda$ at $\lambda=1$ gives \eqref{eq:weak-indices}, while integration gives the converse.
Condition \eqref{eq:flux-indices} implies \eqref{eq:weak-indices} with the same indices $p$ and $q$.
Under \eqref{eq:weak-indices}, however, $\varphi'$ may be discontinuous and $\varphi$ may be linear on some intervals.
Thus, neither classical second differentiability nor strict convexity is required, and minimizers with the same boundary values need not be unique.

Under these weaker assumptions, a natural first goal is H\"older continuity with some exponent determined by the growth indices.
In the scalar variational setting, Lieberman proved H\"older continuity under corresponding growth assumptions; see \cite[Section~6]{Lieberman1991}.
More recently, Behn, Diening, Nowak and Scharle \cite{BDNS2025} proved local boundedness in the vectorial setting and raised the question of H\"older continuity; see \cite[Section~1.1.1, after Theorem~1.1]{BDNS2025}.
Our next theorem answers this question for the autonomous local system considered here and goes beyond a fixed H\"older exponent: we obtain local $W^{1,s}$-regularity for every $s>1$, and hence H\"older continuity with every exponent less than one.

\begin{theorem}\label{thm:main}
Let $\varphi$ be an N-function satisfying
\eqref{eq:potential-scaling}, and let
$u\in W^{1,\varphi}_{\mathrm{loc}}(\Omega;\R^N)$
be a weak solution to \eqref{main_system}.
Then
$u\in W^{1,s}_{\mathrm{loc}}(\Omega;\R^N)$
for every $s\in(1,\infty)$, and hence
$u\in C^{0,\gamma}_{\mathrm{loc}}(\Omega;\R^N)$
for every $\gamma\in(0,1)$.
Moreover, for every $s\in(1,\infty)$, there exists a constant $c>0$ depending only on $n,N,p,q,s$ such that, for every ball $B_{2r}(x_0)\Subset\Omega$,
\begin{equation}\label{eq:reverse-holder}
 \left(\fint_{B_{r/2}(x_0)}|Du|^s\,dx\right)^{1/s}
 \le c\,\varphi^{-1}\left(\fint_{B_{2r}(x_0)}\varphi(|Du|)\,dx\right).
\end{equation}
\end{theorem}

The growth condition \eqref{eq:potential-scaling} can be replaced by the following weaker condition. 
By convexity, this condition implies \eqref{eq:potential-scaling} with possibly different indices.
Since Theorem~\ref{thm:main} holds for every $s\in(1,\infty)$, changing the growth indices does not restrict the range of $s$.

\begin{corollary}\label{rem:weak-growth}
Let $\varphi$ be an N-function satisfying, for some $1<p\le q<\infty$ and $L\ge1$,
\begin{equation}\label{eq:weak-growth}
 L^{-1}\lambda^p\varphi(t)
 \le\varphi(\lambda t)
 \le L\lambda^q\varphi(t)
 \qquad\text{for all }\ \lambda\ge1\quad \text{and}\quad t>0.
\end{equation}
Let $u\in W^{1,\varphi}_{\mathrm{loc}}(\Omega;\R^N)$ be a weak solution to \eqref{main_system}.
Then the conclusions of Theorem~\ref{thm:main} hold, with the constant $c$ depending additionally on $L$.
\end{corollary}

Finally, we develop the argument for Theorem~\ref{thm:main} to obtain local Lipschitz continuity under weaker assumptions than those of Theorem~\ref{thm:gradient}.
We retain \eqref{eq:potential-scaling} and add only a lower growth condition on $\varphi'$.
No second differentiability of $\varphi$ or upper bound on $t\varphi''(t)/\varphi'(t)$ is required; even $\varphi'$ may be discontinuous.
Related Lipschitz estimates were obtained in \cite{Lieberman1993,Marcellini1996,MarcelliniPapi2006,BeckMingione2020,DeFilippisMingione2021} under assumptions controlling the second derivative or the ellipticity ratio.
The finite upper growth index of $\varphi$ is still retained.

\begin{theorem}\label{thm:lower-index-lipschitz}
Let $\varphi$ be an N-function satisfying \eqref{eq:potential-scaling}, and suppose that
\begin{equation}\label{eq:lower-derivative-index}
 \varphi'(\lambda t)\ge\lambda^{p-1}\varphi'(t)
 \qquad\text{for all }\ \lambda\ge1 \quad\text{and}\quad t>0,
 \end{equation}
equivalently,
 $\varphi'(t)/t^{p-1}$ is increasing.
Let $u\in W^{1,\varphi}_{\mathrm{loc}}(\Omega;\R^N)$ be a weak solution to \eqref{main_system}.
Then $u\in W^{1,\infty}_{\mathrm{loc}}(\Omega;\R^N)$, and there exists a constant $c>0$ depending only on $n,N,p,q$ such that, for every ball $B_{2r}(x_0)\Subset\Omega$,
\begin{equation}\label{eq:lower-index-lipschitz}
 \esssup_{B_r(x_0)} |Du|
 \le c\,\varphi^{-1}\left(\fint_{B_{2r}(x_0)}\varphi(|Du|)\,dx\right).
\end{equation}
\end{theorem}

Let us outline the main ideas of the proofs.
For Theorem~\ref{thm:gradient}, Proposition~\ref{prop:baseline} provides the local gradient bound and the supremum decrease estimate, following
\cite[Lemma~5.8 and Theorem~5.9]{DSV2009}.
The main new step is the nondegenerate excess decay in Lemma~\ref{lem:small-excess}.
The classical argument uses a linear comparison system and \eqref{eq:omitted} to estimate the difference between the original nonlinearity and its linear approximation; see \cite[Lemma~6.2]{DSV2009}.
We instead rescale the original nonlinearity directly, without differentiating it, and obtain comparison operators of the form
\[
 \mathcal B(H)=H+[f(\Theta:H)-\Theta:H]\Theta,
 \qquad |\Theta|=1,
\]
where the slopes of $f$ have positive lower and upper bounds depending only on $p$ and $q$.
For a solution $v$ of this system, the scalar function $z=\Theta:Dv$ satisfies
\[
 \diver\big([I_n+(f'(z)-1)\Theta^T\Theta]Dz\big)=0.
\]
This equation is uniformly elliptic even when $f'$ is only measurable.
Scalar H\"older estimates and Schauder estimates for the component equations yield the gradient regularity proved in Lemma~\ref{lem:model}.
Higher integrability controls the difference between the nonlinearities where the normalized gradient is large.
Combining the resulting excess decay with the supremum decrease estimate then proves Theorem~\ref{thm:gradient}, following the scheme of \cite{Uhlenbeck1977,DSV2009}.

For Theorem~\ref{thm:main}, we construct smooth approximations and use Theorem~\ref{thm:gradient} to obtain qualitative regularity of the comparison solutions.
Testing the differentiated system with weighted gradient components is classical; see \cite[Section~1]{Lieberman1993} and \cite[Sections~3--4]{MarcelliniPapi2006}.
The key estimate in Lemma~\ref{lem:angular} is a Caccioppoli inequality for the normalized gradient $Dw/|Dw|$ in which the second derivative of the smooth approximating function does not appear.
Its constant is therefore independent of quantitative bounds on that second derivative. 
Proposition~\ref{prop:moment} then gives uniform local $L^s$ estimates for the gradient for every $s>1$.
The limit and the original solution have the same boundary values and energy, but need not coincide.
Lemma~\ref{lem:flat-interval-transfer} compares their gradient magnitudes and transfers the estimates to the original solution.
A related comparison was used in \cite[Lemma~2.2 and Section~4]{CGGP2017}, where uniform convexity at infinity gives equality of the gradients when the modulus of their sum is sufficiently large.
Our argument instead gives a global comparison of their magnitudes.

For Theorem~\ref{thm:lower-index-lipschitz}, approximation preserves the additional lower growth condition on $\varphi'$. Writing $\Psi$ for the smooth approximating function, the stronger Caccioppoli inequality \eqref{eq:angular-full} then controls derivatives of powers of $\Psi(|Dw|)$. Moser iteration gives a uniform gradient bound without
an upper bound on $t\Psi''(t)/\Psi'(t)$. Strict convexity identifies the limit with the original solution and completes the proof.

Finally, regularity results are also available for suitable classes beyond the growth conditions considered here.
Gradient H\"older regularity for the vectorial energy with density $t\log(1+t)$ was proved in \cite{MingioneSiepe1999}; see also \cite{DeFilippisMingione2023} for scalar nonautonomous problems with nearly linear growth.
Lipschitz estimates covering certain exponential growth conditions were obtained in \cite{Marcellini1996,MarcelliniPapi2006,BeckMingione2020,DeFilippisMingione2021}.
Extending our methods to such cases would require estimates adapted to these growth conditions, possibly together with further approximation arguments.
Here we focus on $(p,q)$-growth with $1<p\le q<\infty$ and do not pursue these extensions.

The paper is organized as follows.
Section~\ref{sec:prelim} contains the preliminaries.
Section~\ref{sec:strong-proof} develops the nonlinear comparison argument and proves Theorem~\ref{thm:gradient}.
Section~\ref{sec:weak-proof} establishes the approximation, Caccioppoli and higher integrability estimates and proves Theorem~\ref{thm:main}.
Section~\ref{sec:lower-index-lipschitz} proves Theorem~\ref{thm:lower-index-lipschitz}.

\section{Preliminaries}\label{sec:prelim}

 For $x_0\in\R^n$ and $r>0$, $B_r(x_0)$ denotes the open ball centered at $x_0$ with radius $r$; when the center is clear, we write $B_r=B_r(x_0)$. For $E\subset\R^n$ with $0<|E|<\infty$, we write
\[
 (f)_E:=\fint_E f\,dx:=\frac{1}{|E|}\int_E f\,dx.
\]
We use $i,j,k\in\{1,\dots,n\}$ for spatial indices and $\alpha,\beta\in\{1,\dots,N\}$ for vector indices, written as subscripts and superscripts, respectively. 
For $P,Q\in\R^{N\times n}$, we use the notation
\[
 P:Q:=\sum_{\alpha=1}^N\sum_{i=1}^n P^\alpha_i Q^\alpha_i,
 \qquad |P|:=\sqrt{P:P}.
\] 
Throughout the paper, $c$ denotes a generic positive constant whose value may vary from line to line.
Subscripts indicate dependence on the corresponding parameters; constants such as $c_n$, $c_b$ and $c_s$ may also vary from line to line.

\begin{definition}\label{def:N}
$\varphi:[0,\infty)\to[0,\infty)$ is called an N-function 
if it is finite and convex, and satisfies
\[
 \varphi(0)=0,\quad \varphi(t)>0\ \ \text{for } \ t>0,
 \quad
 \lim_{t\to0^+}\frac{\varphi(t)}t=0,
 \quad \lim_{t\to\infty}\frac{\varphi(t)}t=\infty.
\]
\end{definition}

The right derivative $\varphi'$ of an N-function is nondecreasing and right-continuous, with $\varphi'(0)=0$. In particular,
\[
 \varphi(t)=\int_0^t\varphi'(s)\,ds.
\]
For convenience, we refer to it simply as the derivative; it agrees with the classical derivative wherever the latter exists. 
Combining \eqref{eq:potential-scaling} with \eqref{eq:weak-indices} yields
\begin{equation}\label{eq:flux-scaling}
 \frac{p}{q}\left(\frac{t}{s}\right)^{p-1}\le\frac{\varphi'(t)}{\varphi'(s)}\le\frac{q}{p}\left(\frac{t}{s}\right)^{q-1}
 \qquad\text{for all }\ 0<s\le t.
\end{equation}
For a bounded open set $D\subset\R^n$, the Orlicz space $L^\varphi(D;\R^N)$ consists of measurable maps $v:D\to\R^N$ such that
\[
 \int_D\varphi\left(\frac{|v|}{\lambda}\right)dx<\infty
 \quad\text{for some }\lambda>0,
\]
and is equipped with the Luxemburg norm
\[
 \|v\|_{L^\varphi(D)}
 :=\inf\left\{\lambda>0:
 \int_D\varphi\left(\frac{|v|}{\lambda}\right)dx\le1\right\}.
\]
The Sobolev--Orlicz space $W^{1,\varphi}(D;\R^N)$ consists of maps $v\in L^\varphi(D;\R^N)$ whose distributional gradient belongs to $L^\varphi(D;\R^{N\times n})$, with norm
\[
 \|v\|_{W^{1,\varphi}(D)} :=\|v\|_{L^\varphi(D)}+\|Dv\|_{L^\varphi(D)}.
\]
We denote by $W^{1,\varphi}_0(D;\R^N)$ the closure of $C_c^\infty(D;\R^N)$ in this norm, and use $W^{1,\varphi}_{\mathrm{loc}}(\Omega;\R^N)$ for the corresponding local space.
Under \eqref{eq:weak-indices}, $\varphi$ and its conjugate $\varphi^*(a):=\sup_{t\ge0}\{at-\varphi(t)\}$ satisfy the $\Delta_2$-condition.
Consequently, these spaces are reflexive, and $v\in L^\varphi(D)$ if and only if $\int_D\varphi(|v|)\,dx<\infty$.
For these standard facts, we refer to \cite{HarjulehtoHasto2019}.

For an N-function $\varphi$, define $A:\R^{N\times n}\to\R^{N\times n}$ by
\[
 A(P):=\varphi'(|P|)\frac{P}{|P|}.
 \]
Under \eqref{eq:weak-indices}, convexity and the weak formulation imply that every weak solution $u$ to \eqref{main_system} minimizes the energy among maps with the same boundary values:
\[
 \int_{B_r}\varphi(|Du|)\,dx
 \le\int_{B_r}\varphi(|Dv|)\,dx
 \qquad
 \text{for every }\ v\in u+W^{1,\varphi}_0(B_r;\R^N),
\]
whenever $B_r\Subset\Omega$.

We also record the following standard estimates; see \cite{DSV2009}. They remain valid under \eqref{eq:flux-indices}, since $\varphi'$ is locally absolutely continuous and satisfies \eqref{eq:indices} almost everywhere.

\begin{lemma}\label{lem:secant}
Let $\varphi$ be an N-function satisfying \eqref{eq:flux-indices}, and let $V$ be defined in \eqref{eq:AV}.
There exists a constant $c\ge1$, depending only on
$n,N,p,q$, such that for all
$P,Q\in\R^{N\times n}$ with $|P|+|Q|>0$,
\[
 (A(P)-A(Q)):(P-Q) 
 \ge c^{-1} \frac{\varphi'(|P|+|Q|)}{|P|+|Q|}|P-Q|^2,
\]
\[
 c^{-1}\frac{\varphi'(|P|+|Q|)}{|P|+|Q|}|P-Q|^2
 \le |V(P)-V(Q)|^2
 \le c\frac{\varphi'(|P|+|Q|)}{|P|+|Q|}|P-Q|^2,
\]
and
\[
 |A(P)-A(Q)| \le c\frac{\varphi'(|P|+|Q|)}{|P|+|Q|}|P-Q|.
\]
\end{lemma}

\section{Gradient H\"older regularity under strong growth}\label{sec:strong-proof}

We first record uniform estimates for the system
\eqref{main_system} with an N-function $\varphi$ satisfying
\eqref{eq:flux-indices}. The local gradient bound and the
supremum decrease estimate were proved in
\cite[Lemma~5.8 and Theorem~5.9]{DSV2009} under the additional
assumption that $\varphi\in C^2((0,\infty))$. We explain below
how to adapt the proof under \eqref{eq:flux-indices}, using
the shifted problems from \cite[Section~4]{DSV2009}.

\begin{proposition}\label{prop:baseline}
Let $\varphi\in C^1([0,\infty))$ be an N-function satisfying
\eqref{eq:flux-indices}, and let
$u\in W^{1,\varphi}_{\mathrm{loc}}(\Omega;\R^N)$ be a weak
solution to \eqref{main_system}. There exists a constant
$c\ge1$, depending only on $n,N,p$ and $q$, such that,
for every ball $B_{2r}(x_0)\Subset\Omega$,
\begin{equation}\label{eq:Lipschitz}
\begin{aligned}  
\esssup_{B_r(x_0)} \varphi(|Du|)
& \le c\,\fint_{B_{2r}(x_0)}\varphi(|Du|)\,dx\\
& \le c\left(\mathcal E(u,B_{2r}(x_0)) +\big|(V(Du))_{B_{2r}(x_0)}\big|\right),
\end{aligned}\end{equation}
and
\begin{equation}\label{eq:sup-drop}
 \mathcal E(u,B_{r/2}(x_0))
 \le c\bigg(\esssup_{B_{r}(x_0)} \varphi(|Du|)-\esssup_{B_{r/2}(x_0)} \varphi(|Du|) \bigg).
\end{equation}
\end{proposition}

\begin{proof}[Sketch of proof]
We follow \cite[Sections~4--5]{DSV2009}, explaining the changes
needed when $\varphi'$ is only locally absolutely continuous.
Throughout the proof, all balls are centered at $x_0$, and
we omit the center from the notation. Choose $\rho>2r$ such
that $B_\rho\Subset\Omega$. For $\lambda>0$, set
\[
 \varphi_\lambda(t)
 :=\int_0^t\frac{\varphi'(\lambda+s)}{\lambda+s}s\,ds,
 \qquad
 A_\lambda(P)
 :=\frac{\varphi'(\lambda+|P|)}{\lambda+|P|}P,
\]
and let $u_\lambda$ be the minimizer with boundary values $u$
in $B_\rho$. Note that $u_\lambda$ is a weak solution to
$\diver A_\lambda(Du_\lambda)=0$ in $B_\rho$.
By \eqref{eq:flux-indices},
\[
 h_\lambda(t)
 :=\frac{t\varphi_\lambda''(t)}{\varphi_\lambda'(t)}
 =\frac{t}{\lambda+t} \frac{(\lambda+t)\varphi''(\lambda+t)}{\varphi'(\lambda+t)}+\frac{\lambda}{\lambda+t}
\]
satisfies
\[
 \min\{1,p-1\}\le h_\lambda(t)\le\max\{1,q-1\} \qquad\text{for a.e. }t>0,
\]
uniformly in $\lambda$. Here we choose a representative of $h_\lambda$ satisfying these inequalities for every $t>0$.

For each fixed $\lambda>0$, $\varphi_\lambda\in W^{2,\infty}_{\mathrm{loc}}([0,\infty))$, and $A_\lambda$ is uniformly elliptic on bounded sets, with ellipticity constants that may depend on $\lambda$ and the set.
As in \cite[Proposition~5.6]{DSV2009}, \cite[Theorem~A]{MarcelliniPapi2006} gives $Du_\lambda\in L^\infty_{\mathrm{loc}}(B_\rho)$.
Difference quotients then yield $u_\lambda\in W^{2,2}_{\mathrm{loc}}(B_\rho;\R^N)$.

The first change is in \cite[Lemma~5.3]{DSV2009}:
we use the Sobolev chain rule for the locally Lipschitz map
$A_\lambda$, instead of the chain rule for a $C^1$ map.
The resulting identities hold almost everywhere.
Set
\[
 w_\lambda:=\varphi_\lambda(|Du_\lambda|)
 \qquad\text{and}\qquad
 V_\lambda(P) :=\sqrt{\frac{\varphi_\lambda'(|P|)}{|P|}}P.
\]
The argument of \cite[Lemma~5.4]{DSV2009} then gives
\cite[(5.7)]{DSV2009}, namely,
\[
 -\int_{B_\rho}G_\lambda(Du_\lambda)Dw_\lambda
       \cdot D\eta\,dx
 \ge c\int_{B_\rho}
       \eta\,|D(V_\lambda(Du_\lambda))|^2\,dx,
\]
where $0\le\eta\in C_c^\infty(B_\rho)$ and 
\[
 G_\lambda(P) :=I+(h_\lambda(|P|)-1)\frac{P^TP}{|P|^2}.
\]
These matrices are uniformly elliptic with constants depending only on
$p$ and $q$. Following the proofs of
\cite[Lemma~5.8 and Theorem~5.9]{DSV2009}, we obtain
the local gradient bound and the supremum decrease estimate
for $u_\lambda$, with constants independent of $\lambda$.

The second change concerns the limit in \cite[Theorem~5.5]{DSV2009}, where the continuity of $\varphi''$ is used to identify the limiting coefficient matrix.
Under \eqref{eq:flux-indices}, $\varphi''$ need not be continuous, so this passage is not directly available.
Instead, we first derive the estimates for $u_\lambda$ and then pass to the limit in these estimates, without requiring convergence of the coefficient matrices.
The energy and compactness arguments of \cite[Sections~3--4]{DSV2009} give, along a subsequence,
\[
 Du_\lambda\to Du\quad\text{a.e. in }B_\rho,
 \qquad\text{and}\qquad
 V_\lambda(Du_\lambda)\to V(Du) \quad\text{in }L^2_{\mathrm{loc}}(B_\rho).
\]
All subsequent limits are taken along this subsequence.
The uniform local gradient bound consequently gives $w_\lambda\to w:=\varphi(|Du|)$ in $L^s_{\mathrm{loc}}(B_\rho)$ for every finite $s\ge1$.
As in the proof of \cite[Lemma~5.8, (5.15)]{DSV2009},
\[
 \esssup_{B_r}w_\lambda \le c\fint_{B_{2r}}w_\lambda\,dx,
\]
where $c$ is independent of $\lambda$.
Passing to the limit proves the first inequality in \eqref{eq:Lipschitz}. The second follows from $\varphi(|P|)\le c|V(P)|^2$ and the definition of $\mathcal E$.

To prove \eqref{eq:sup-drop}, set $M_r:=\esssup_{B_r}w$. For $0<s<r$, the local boundedness estimate for the subsolution $(w_\lambda-M_r)_+$ gives
\[
 \esssup_{B_s}(w_\lambda-M_r)_+
 \le C_{s,r}\fint_{B_r}(w_\lambda-M_r)_+\,dx
 \longrightarrow0
 \quad\text{as }\lambda\to0.
\]
Hence, using also the almost everywhere convergence,
\[
 \limsup_{\lambda\to0}\esssup_{B_s}w_\lambda\le M_r
 \quad\text{and}\quad
 \esssup_{B_{s/2}}w \le\liminf_{\lambda\to0}\esssup_{B_{s/2}}w_\lambda.
\]
Writing $\mathcal E_\lambda$ for the excess defined using $V_\lambda$, we pass to the limit in
\[
 \mathcal E_\lambda(u_\lambda,B_{s/2})
 \le c\left(\esssup_{B_s}w_\lambda-\esssup_{B_{s/2}}w_\lambda \right)
\]
and obtain
\[
 \mathcal E(u,B_{s/2}) \le c\left(M_r-\esssup_{B_{s/2}}w\right).
\]
Letting $s\nearrow r$, we obtain \eqref{eq:sup-drop},
since the excess and the essential supremum converge to their corresponding values on $B_{r/2}$.
\end{proof}

We establish $C^{1,\gamma_0}$-regularity for an elliptic system
of a special form, which will be used in the proof of
Theorem~\ref{thm:gradient}. While the usual approach uses
linearization to obtain a linear limiting system, our proof
retains a nonlinear system in the limit. The key observation
is that a certain linear combination of the gradient components
satisfies a scalar uniformly elliptic equation.

\begin{lemma}\label{lem:model}
Let $\Theta\in\R^{N\times n}$ with $|\Theta|=1$, and let $f:\R\to\R$ satisfy $f(0)=0$ and
\begin{equation}\label{eq:f-slope}
 p-1\le\frac{f(t)-f(s)}{t-s}\le q-1
 \qquad\text{for all }\ s,t\in\R\ \text{ with }\ s<t.
\end{equation}
Define
\[
 \mathcal B(P):=P+[f(\Theta:P)-\Theta:P]\Theta.
\]
Then every weak solution $v\in W^{1,2}(B_1;\R^N)$ to
\[
 \diver\mathcal B(Dv)=0\quad\text{in }B_1
\]
belongs to $C^{1,\gamma_0}(B_{1/2};\R^N)$ for some $\gamma_0\in(0,1)$ depending only on $n,N,p,q$, and there exists a constant $c>0$, with the same dependence, such that
\[
 \sup_{B_{1/2}}|Dv|+[Dv]_{C^{0,\gamma_0}(B_{1/2})}
 \le c\left(\fint_{B_1}|Dv|^2\,dx\right)^{1/2}.
\]
\end{lemma}

\begin{proof}
Let $P,Q\in \R^{N\times n}$.
Note that
\[
 \mathcal B(P)-\mathcal B(Q)
 =\underbrace{P-Q-[\Theta:(P-Q)]\Theta}_{=:Z}
  +\big(f(\Theta:P)-f(\Theta:Q)\big)\Theta.
\]
Since $|\Theta|=1$, we have $\Theta:Z=0$ and
$|P-Q|^2=[\Theta:(P-Q)]^2+|Z|^2$.
By \eqref{eq:f-slope}, we obtain
\begin{equation}\label{eq:B_mono}\begin{aligned}
 (\mathcal B(P)-\mathcal B(Q)):(P-Q)
 &=|Z|^2+\big(f(\Theta:P)-f(\Theta:Q)\big)
       [\Theta:(P-Q)]\\
 &\ge |Z|^2+(p-1)[\Theta:(P-Q)]^2\\
 &\ge \min\{1,p-1\}|P-Q|^2,
\end{aligned}\end{equation}
and
\[\begin{aligned}
 |\mathcal B(P)-\mathcal B(Q)|^2
 &=|Z|^2+\big|f(\Theta:P)-f(\Theta:Q)\big|^2\\
 &\le |Z|^2+(q-1)^2[\Theta:(P-Q)]^2\\
 &\le\max\{1,q-1\}^2|P-Q|^2.
\end{aligned}\]
Therefore, a standard difference-quotient argument yields
$v\in W^{2,2}_{\mathrm{loc}}(B_1;\R^N)$.

Set
$
 z:=\Theta:Dv
 =\sum_{\alpha=1}^N\sum_{i=1}^n
   \Theta^\alpha_i\partial_i v^\alpha
 \in W^{1,2}_{\mathrm{loc}}(B_1).
$
Since
\begin{equation}\label{eq:B1}
 \diver\big([\mathcal B(Dv)]^\alpha\big)
 =\diver\big(Dv^\alpha+(f(z)-z)\Theta^\alpha\big)=0
 \quad\text{for each }\alpha=1,\dots,N,
\end{equation}
where $[\mathcal B(Dv)]^\alpha$ and
$\Theta^\alpha=[\Theta^\alpha_1\,\cdots\,\Theta^\alpha_n]$
are the $\alpha$-th rows of $\mathcal B(Dv)$ and $\Theta$,
respectively, applying
$\sum_{k=1}^n\Theta^\alpha_k\partial_k$
and summing over $\alpha$ yield
\begin{align*}
 0
 &=\diver\bigg( D\Big[\underbrace{\sum_{\alpha,k} \Theta^\alpha_k\partial_kv^\alpha}_{=z}\Big] +(f'(z)-1) \underbrace{\sum_{\alpha,k}c\Theta^\alpha_k\partial_kz\,\Theta^\alpha}_{=\Theta^T\Theta Dz}\bigg)\\
 &=\diver\big([I_n+(f'(z)-1)\Theta^T\Theta]Dz\big)
\end{align*}
in the sense of distributions in $B_1$. Note that, from \eqref{eq:f-slope}, $f$ is differentiable a.e. with $p-1\le f'(t)\le q-1$. 
Since $|\Theta|=1$, for every $\xi\in\R^n$ we have
\[
 \min\{1,p-1\}|\xi|^2
 \le [I_n+(f'(z)-1)\Theta^T\Theta]\xi\cdot\xi
 \le\max\{1,q-1\}|\xi|^2.
\]
The De Giorgi--Nash--Moser estimate
\cite[Theorem~8.22]{GT} gives
\[
 \sup_{B_{3/4}}|z| +[z]_{C^{0,\gamma_0}(B_{3/4})}
 \le c\|z\|_{L^2(B_{7/8})}
 \le c\|Dv\|_{L^2(B_1)}.
\]
Since $f(0)=0$ and $f$ is Lipschitz, we obtain
\[
 \sup_{B_{3/4}}|f(z)-z| +[f(z)-z]_{C^{0,\gamma_0}(B_{3/4})}
 \le c\|Dv\|_{L^2(B_1)}.
\]
Finally, applying the interior gradient Schauder estimate for the Poisson equation with H\"older continuous divergence data \cite[Theorem~8.32]{GT} to
\[
 \Delta (v^\alpha-(v^\alpha)_{B_{3/4}}) = \Delta v^\alpha
 =-\diver\big([f(z)-z]\Theta^\alpha\big),
 \qquad \alpha=1,\dots,N,
\]
which follows from \eqref{eq:B1}, we obtain
\[\begin{aligned}
\sup_{B_{1/2}}|Dv|+[Dv]_{C^{0,\gamma_0}(B_{1/2})}
&\le c\left(\|v-(v)_{B_{3/4}}\|_{L^2(B_{3/4})}+\|f(z)-z\|_{C^{0,\gamma_0}(B_{3/4})}\right)\\
&\le c\|Dv\|_{L^2(B_1)},
\end{aligned}\]
where the last inequality follows from Poincar\'e's inequality and the preceding estimate.
\end{proof}

\begin{lemma}\label{lem:small-excess}
Let $\varphi$ be an N-function satisfying
\eqref{eq:flux-indices}, and let
$u\in W^{1,\varphi}_{\mathrm{loc}}(\Omega;\R^N)$ be a weak
solution to \eqref{main_system}. There exist constants
$c_*\ge1$ and $\gamma_0\in(0,1)$, depending only on $n,N,p$
and $q$, such that for every $\tau\in(0,1/16)$ there exists
$\epsilon_*(\tau)>0$, depending only on $n,N,p,q$ and $\tau$,
for which, on every ball $B_r\Subset\Omega$,
\[
 \mathcal E(u,B_r)\le\epsilon_*(\tau)|(V(Du))_{B_r}|^2
 \quad\Longrightarrow\quad
 \mathcal E(u,B_{\tau r})
 \le c_*\tau^{2\gamma_0}\mathcal E(u,B_r).
\]
Here $B_r$ and $B_{\tau r}$ have the same center.
\end{lemma}

\begin{proof}
\noindent\textit{Step 1. Normalization and higher integrability.}
By translating and rescaling, we may
assume that $B_r=B_1=B_1(0)$. The case $\mathcal E(u,B_1)=0$ is immediate, since $V(Du)$ is constant almost everywhere in $B_1$, and hence $\mathcal E(u,B_\tau)=0$.  Assume that
$\mathcal E(u,B_1)>0$ and that the smallness condition $\mathcal E(u,B_1)\le\epsilon_*(\tau)|(V(Du))_{B_1}|^2$ holds,
where $\epsilon_*(\tau)$ will be chosen below.

Define
\begin{gather*}
 Q:=V^{-1}((V(Du))_{B_1}),
 \quad U:=\frac u{|Q|},
 \quad \Theta:=\frac Q{|Q|},\\
 \Phi(t):=\frac{\varphi(|Q|t)}{|Q|\varphi'(|Q|)},
 \quad\text{and}\quad
 \epsilon^2:=
 \frac{\mathcal E(u,B_1)}{|(V(Du))_{B_1}|^2}
 =\frac{\mathcal E(u,B_1)}{|Q|\varphi'(|Q|)}\le \epsilon_*(\tau).
\end{gather*}
Then $\Phi$ satisfies \eqref{eq:flux-indices} with the
same indices $p$ and $q$, $\Phi'(t)=\phi'(|Q|t)/\phi'(|Q|)$, and hence $\Phi'(1)=1$. Set
\[
 \mathcal A(P):=\Phi'(|P|)\frac P{|P|}
 \quad\text{and}\quad
 \mathcal V(P):=\sqrt{\frac{\Phi'(|P|)}{|P|}}P.
\]
Direct substitution gives
\begin{equation}\label{lem:small-excess_pf0}
 \diver\mathcal A(DU)= \diver\left(\frac{\phi'(|Du|)}{\phi'(|Q|)} \frac{Du}{|Du|}\right)= \frac{1}{\phi'(|Q|)} \diver\left(\frac{\phi'(|Du|)}{|Du|} Du\right) =0.
\end{equation}
Furthermore,
\[
\mathcal V(DU)=\frac{V(Du)}{\sqrt{|Q|\varphi'(|Q|)}} 
\quad\text{and}\quad   (\mathcal V(DU))_{B_1}=\Theta =\mathcal{A}(\Theta),
 \]
which implies
 \[
 \fint_{B_1}|\mathcal V(DU)-  \Theta|^2\,dx 
 =\fint_{B_1}|\mathcal V(DU)- (\mathcal V(DU))_{B_1} |^2\,dx =\frac{\mathcal E (u,B_1)}{|Q|\phi'(|Q|)}
 =\epsilon^2.
\]
We choose $\epsilon_*(\tau)\le1$, so that $\epsilon\le1$.
Consequently, applying the local gradient bound \eqref{eq:Lipschitz} in Proposition~\ref{prop:baseline} to the system \eqref{lem:small-excess_pf0} and the previous identities with $|\Theta|=1$, 
\begin{equation}\label{lem:small-excess_pf1}
 \esssup_{B_{1/2}} |DU| \le c\, \Phi^{-1}\left(\fint_{B_1}|\mathcal V(DU)|^2\,dx\right) \le c_L,
\end{equation}
where $c_L\ge1$ depends only on $n,N,p$ and $q$.

We next define 
\[
 v(x):=\frac{U(x)-\Theta x}{\epsilon},
 \quad 
 \mathcal B_\epsilon(P)
 :=\frac{\mathcal A(\Theta+\epsilon P)-\Theta}{\epsilon}=\frac{\mathcal A(\Theta+\epsilon P)-\mathcal A(\Theta)}{\epsilon}.
\]
Then, by Lemma~\ref{lem:secant} and \eqref{lem:small-excess_pf1},  
\[\begin{aligned}
 \int_{B_{1/2}}|Dv|^2\,dx 
 &= \frac{1}{\epsilon^2} \int_{B_{1/2}}|DU-\Theta|^2\,dx
 \le\frac{c}{\epsilon^2} \int_{B_{1/2}}|\mathcal V(DU)-\Theta|^2\,dx \\
 &\le\frac{c}{\epsilon^2} \int_{B_{1/2}}\frac{|DU|+1}{\Phi'(|DU|+1)} |\mathcal V(DU)-\Theta|^2\,dx\\
&\le\frac{c}{\epsilon^2} \int_{B_{1/2}} |\mathcal V(DU)-\Theta|^2\,dx \le c.
\end{aligned}\]
Moreover, by \eqref{lem:small-excess_pf0},
\begin{equation}\label{lem:small-excess_pf2}
 \diver\mathcal B_\epsilon(Dv)= \diver\left(\frac{\mathcal A(DU)-\Theta}{\epsilon}\right) =0\quad\text{in }B_{1/2},
{\color{red} \qquad
} \end{equation}
and Lemma~\ref{lem:secant}, applied to $\Phi$ in place of  $\phi$ with $P=DU$ and $Q=\Theta$, together with \eqref{lem:small-excess_pf1} gives
\[\begin{aligned}
 \mathcal B_\epsilon(Dv):Dv= \frac1{\epsilon^2} (\mathcal A(DU)-\mathcal A(\Theta)):(DU-\Theta)
 &\ge \frac{1}{c}\frac{\Phi'(|DU|+1)}{|DU|+1}|Dv|^2\\
 &\ge \frac{1}{c}|Dv|^2
 \end{aligned}\]
 and
\[
|\mathcal B_\epsilon(Dv)|=\frac{|\mathcal A(DU) - \mathcal A (\Theta)|}{\epsilon} \le c \frac{\Phi'(|DU|+1)}{|DU|+1} |Dv| \le c|Dv|.
\]
Therefore, by testing the system \eqref{lem:small-excess_pf2} with the test function $\eta^2 (v-(v)_{B_{2\rho}(y)}) $, where $B_{2\rho}(y)\subset B_{1/2}$ and  $\eta\in C^\infty_{0}(B_{2\rho}(y))$ with $0\le\eta\le 1$, $\eta\equiv 1$ on $B_{\rho}(y)$, and $|D\eta|\le 2/\rho$, and using the previous inequalities, we obtain that
\[
 \fint_{B_\rho(y)}|Dv|^2\,dx
 \le\frac c{\rho^2}
       \fint_{B_{2\rho}(y)}|v-(v)_{B_{2\rho}(y)}|^2\,dx \le c \left(\fint_{B_{2\rho}(y)}|Dv|^{2_*}\,dx\right)^{\frac{2}{2_*}}.
\]
Here $2_*:=\max\{2n/(n+2),3/2\}\in(1,2)$, and the last inequality follows from the Sobolev--Poincar\'e inequality.
Therefore, by Gehring's lemma \cite{Gehring1973,Giaquinta},
there exists $\delta_0>0$, depending only on $n,N,p$ and $q$,
such that
\begin{equation}\label{eq:Dv_high}
 \left(\fint_{B_{3/8}}|Dv|^{2+\delta_0}\,dx\right)^{2/(2+\delta_0)}\le c\fint_{B_{1/2}}|Dv|^{2}\,dx \le c.
\end{equation}

\medskip
\noindent\textit{Step 2. Construction of the comparison system.}
We next construct a nonlinear system to which
Lemma~\ref{lem:model} applies and compare its solution with $v$.
Set
\begin{equation}\label{eq:m}
 m:=\epsilon^{-2/(4+\delta_0)},
 \end{equation}
 where $\delta_0>0$ is given in \eqref{eq:Dv_high}.
Since $\Phi$ satisfies \eqref{eq:flux-indices} and $\Phi'(1)=1$, we have  
\begin{equation}\label{lem:small-excess_pf3}
 \frac12\le\frac{\Phi'(1+t)}{1+t}\le2
 \qquad\text{whenever }\ |t|\le 2\epsilon m = 2\epsilon^{\frac{2+\delta_0}{4+\delta_0}}<\frac12,
\end{equation}
provided that $\epsilon$ is sufficiently small depending only on
$n,N,p,q$. 

Define
\[
 f_\epsilon(s):= 
\begin{cases} 
(p-1) \left(s+2m\right)+\frac{\Phi'(1-2\epsilon m)-1}{\epsilon} &\text{if }\ s< -2m \\
 \vspace{-0.2cm}\ & \ \\
\displaystyle \ \ \ \ \ \ \ \ \ \ \  \frac{\Phi'(1+\epsilon s)-1}{\epsilon} & \text{if }\  -2m  \le s \le 2m \\
 \vspace{-0.2cm} \ & \ \\
(p-1) \left(s-2m\right)+\frac{\Phi'(1+2\epsilon m)-1}{\epsilon} &\text{if }\ s> 2m .
\end{cases}
\]
Since $\Phi'$ is locally absolutely continuous and satisfies \eqref{eq:flux-indices}, $f_\epsilon$ is  absolutely continuous and, by \eqref{lem:small-excess_pf3}, satisfies
\[
 \frac{p-1}{2}
 \le(p-1)\frac{\Phi'(1+\epsilon s)}{1+\epsilon s}
 \le f_\epsilon'(s)
 \le(q-1)\frac{\Phi'(1+\epsilon s)}{1+\epsilon s}
 \le2(q-1)
\]
for almost every $s\in[-2m,2m]$. Moreover, $f_\epsilon'(s)=p-1$ for all $s\in \R$ with $|s|>2m$. Therefore, we obtain 
\[
 \frac{p+1}{2} -1 =\frac{p-1}{2}
 \le \frac{f_\epsilon(t)-f_\epsilon(s)}{t-s}
 \le 2(q-1)= (2q-1) -1
\]
for all $s,t\in\R$ with $s<t$. Define
\[
 \widetilde{\mathcal B}_\epsilon(P)
 :=P+[f_\epsilon(\Theta:P)-\Theta:P]\Theta.
\]

We estimate the difference between $\mathcal B_\epsilon$
and $\widetilde{\mathcal B}_\epsilon$. Suppose that $|P| \le m$ in this paragraph. Since $ |\Theta+\epsilon P|^2= 1+ 2\epsilon (\Theta:P) + \epsilon^2|P|^2 $, Taylor's formula gives 
\begin{align*}
 |\Theta+\epsilon P| =1+\epsilon(\Theta:P)+R_1
\ \ \text{and}\ \ 
 \frac{\Theta+\epsilon P}{|\Theta+\epsilon P|} =\Theta+\epsilon\big(P-(\Theta:P)\Theta\big)+R_2,
\end{align*}
where
$|R_1|+|R_2|\le c\epsilon^2|P|^2$ with 
the constant $c>0$ depending only on $n$ and $N$,
since $\epsilon|P|\le \epsilon m \le1/4$.
Then 
\[\begin{aligned}
 \Phi'(|\Theta+\epsilon P|)
 &=1+\epsilon f_\epsilon\left(\frac{|\Theta+\epsilon P|-1}{\epsilon}\right)
 =1+\epsilon f_\epsilon(\Theta:P+R_1/\epsilon)\\
 &= 1+\epsilon f_\epsilon(\Theta:P)+\underbrace{\epsilon\left(f_\epsilon(\Theta:P+R_1/\epsilon)- f_\epsilon(\Theta:P)\right)}_{=:R_3},
\end{aligned}\]
where $|R_3|\le 2 (q-1)|R_1|\le c \epsilon^2|P|^2$. Multiplying by the expansion of
$(\Theta+\epsilon P)/|\Theta+\epsilon P|$ and using
$|f_\epsilon(\Theta:P)|\le2(q-1)|P|$, we obtain
\[\begin{aligned}
 \mathcal A(\Theta+\epsilon P)
 &=  \Phi'(|\Theta+\epsilon P|) \frac{\Theta+\epsilon P}{|\Theta+\epsilon P|} \\
 &= (1+\epsilon f_\epsilon(\Theta:P)+R_3)(\Theta+\epsilon\big(P-(\Theta:P)\Theta\big)+R_2)\\
 &=\Theta+\epsilon(P+[f_\epsilon(\Theta:P)-\Theta:P]\Theta)+R_4\\
 &=\Theta+\epsilon\widetilde{\mathcal B}_\epsilon(P)+R_4,
\end{aligned}\]
where $|R_4|\le c\epsilon^2|P|^2$.
Consequently,
\[
 |\mathcal B_\epsilon(P)
       -\widetilde{\mathcal B}_\epsilon(P)|
 \le c\epsilon|P|^2
 \qquad\text{for }\ |P|\le m.
\]
Using this estimate on $\{|Dv|\le m\}$, \eqref{eq:Dv_high} and \eqref{eq:m}, we conclude that
\begin{equation}\label{eq:comparison_BDv}
\begin{aligned}
 \int_{B_{1/4}} |\mathcal B_\epsilon(Dv)-\widetilde{\mathcal B}_\epsilon(Dv)|^2\,dx
 &\le c\epsilon^2m^4 +c\int_{B_{1/4}\cap\{|Dv|>m\}}|Dv|^2\,dx \\
 &\le c\big(\epsilon^2m^4+m^{-\delta_0}\big) \le c\epsilon^\nu,
 \end{aligned}
\end{equation}
where $\nu:=\frac{2\delta_0}{4+\delta_0}$.

\medskip
\noindent\textit{Step 3. Comparison and excess decay.} Let $h\in v+W^{1,2}_0(B_{1/4};\R^N)$ be the unique weak
solution to
\[
 \diver\widetilde{\mathcal B}_\epsilon(Dh)=0
 \qquad\text{in }\ B_{1/4}.
\]
Note that the existence and uniqueness of $h$ follow from the strong monotonicity
and Lipschitz continuity of $\widetilde{\mathcal B}_\epsilon$, as in the proof of Lemma~\ref{lem:model}.
Testing the system by $h-v$ gives
\[
 \int_{B_{1/4}}|Dh|^2\,dx
 \le c\int_{B_{1/4}}|Dv|^2\,dx\le c.
\]
Applying Lemma~\ref{lem:model}, after rescaling, to $\mathcal B=\widetilde{\mathcal B}_\epsilon$ and $f=f_\epsilon$ with $(p,q)$ replaced by $(\frac{p+1}{2},2q-1)$, we obtain $\gamma_0\in(0,1)$ and $c\ge1$, depending only on
$n,N,p$ and $q$, such that
\[
 \sup_{B_{1/8}}|Dh| + [Dh]_{C^{0,\gamma_0}(B_{1/8})}\le c.
\]
Subtracting the systems for $v$ and $h$, and testing by
$v-h$, we obtain
\begin{align*}
\int_{B_{1/4}} &[\widetilde{\mathcal B}_\epsilon(Dv)-\widetilde{\mathcal B}_\epsilon(Dh)]:(Dv-Dh)\,dx\\
& = \int_{B_{1/4}} [\widetilde{\mathcal B}_\epsilon(Dv)-\mathcal B_\epsilon(Dv)] :(Dv-Dh)\,dx.
\end{align*}
By \eqref{eq:B_mono}, Young's inequality, and \eqref{eq:comparison_BDv}, we obtain
\begin{equation}\label{eq:comparison_DhDv}
 \int_{B_{1/4}}|Dv-Dh|^2\,dx
 \le c\int_{B_{1/4}}
 |\widetilde{\mathcal B}_\epsilon(Dv)
       -\mathcal B_\epsilon(Dv)|^2\,dx
 \le c\epsilon^\nu.
\end{equation}

We finally prove the excess decay estimate. For $\tau\in(0,1/16)$, it follows from the preceding H\"older estimate and \eqref{eq:comparison_DhDv} that
\begin{align*}
 \fint_{B_\tau}|Dv-(Dh)_{B_\tau}|^2\,dx
 &\le c\fint_{B_\tau}|Dv-Dh|^2\,dx +c\fint_{B_\tau}|Dh-(Dh)_{B_\tau}|^2\,dx\\
& \le c\big(\tau^{-n}\epsilon^\nu+\tau^{2\gamma_0}\big).
\end{align*}

The preceding $L^\infty$-estimate for $Dh$ and $|\Theta|=1$ also give
\[
 \frac12\le 1-\epsilon|(Dh)_{B_\tau}| \le|\Theta+\epsilon(Dh)_{B_\tau}| \le 1+\epsilon|(Dh)_{B_\tau}| \le2,
\]
by assuming $\epsilon$ sufficiently small depending only on $n,N,p,q$.
Together with \eqref{lem:small-excess_pf1}, Lemma~\ref{lem:secant} implies
\[\begin{aligned}
 |\mathcal V(DU)-&\mathcal V(\Theta+\epsilon(Dh)_{B_\tau})|^2 \\
       &\le c \frac{\Phi'(|DU|+|\Theta+\epsilon(Dh)_{B_\tau}|)}{|DU|+|\Theta+\epsilon(Dh)_{B_\tau}|} |DU-\Theta-\epsilon(Dh)_{B_\tau}|^2\\
 &\le c\epsilon^2|Dv-(Dh)_{B_\tau}|^2.
\end{aligned}\]
Using this and recalling the definition of $\epsilon$ and $\mathcal V$ at the beginning of the proof, we obtain
\begin{align*}
 \frac{\mathcal E(u,B_\tau)}{\mathcal E(u,B_1)} 
 &=\frac{|Q|\phi'(|Q|)}{\mathcal E(u,B_1)} \fint_{B_\tau} |\mathcal V(DU)-(\mathcal V(DU))_{B_\tau}|^2\,dx\\
 &=\frac1{\epsilon^2}\fint_{B_\tau}|\mathcal V(DU)-(\mathcal V(DU))_{B_\tau}|^2\,dx\\
 &\le\frac1{\epsilon^2}\fint_{B_\tau}|\mathcal V(DU)-\mathcal V(\Theta+\epsilon(Dh)_{B_\tau})|^2\,dx\\
 &\le c\fint_{B_\tau}|Dv-(Dh)_{B_\tau}|^2\,dx\\
 &\le c\big(\tau^{-n}\epsilon^\nu+\tau^{2\gamma_0}\big).
\end{align*}

At this stage we choose 
\[
 \epsilon_*(\tau):= \min\{     \epsilon_0,     \tau^{2(n+2\gamma_0)/\nu}\},
\]
where $\epsilon_0\in (0,1]$ is a sufficiently small constant depending only on $n,N,p,q$, chosen so that all the preceding estimates hold under the condition that $\epsilon^2\le \epsilon_0$. 
The condition gives
$\tau^{-n}\epsilon^\nu\le\tau^{-n}\epsilon_*(\tau)^{\nu/2} \le\tau^{2\gamma_0}$.
Therefore,
\[
 \mathcal E(u,B_\tau)
 \le c_* \tau^{2\gamma_0}\mathcal E(u,B_1).
\]
Rescaling back to $B_r$ completes the proof.
\end{proof}

We now prove the H\"older continuity of the gradient by
combining the small-excess estimate with the supremum
decrease estimate and following the iteration argument
in \cite[Section~6]{DSV2009}.

\begin{proof}[Proof of Theorem~\ref{thm:gradient}]
Proposition~\ref{prop:baseline} and Lemma~\ref{lem:small-excess} apply directly to $u$.
Let $B_r\Subset\Omega$ and write
\[
 M(\rho):=\esssup_{B_\rho}\varphi(|Du|), 
 \qquad \rho\in(0,r].
\]
Fix $\tau\in(0,1/16)$ so that
\[
 c_*\tau^{2\gamma_0}\le\frac14.
\]
Let $c_1\ge1$ be a common constant in
Proposition~\ref{prop:baseline}, and set
\[
 \epsilon_0:=
 \min\left\{1,\frac1{2c_1},
                  \frac{\epsilon_*(\tau)}{2c_1}\right\}
 \qquad\text{and}\qquad
 \delta:=1-\frac{\epsilon_0(4\tau)^n}{8c_1}
 \in(1/2,1).
\]

We first verify the two alternatives. If
$\mathcal E(u,B_r) \le \epsilon_0\esssup_{B_{r/2}} \varphi(|Du|)$, then \eqref{eq:Lipschitz} and the choice of $\epsilon_0$ imply
\[\begin{aligned}
 \mathcal E(u,B_r) \le \epsilon_0M(r/2) 
&\le c_1\epsilon_0 \left(\mathcal E(u,B_r)+|(V(Du))_{B_r}|^2\right) \\
&\le \frac12\mathcal E(u,B_r)+c_1\epsilon_0|(V(Du))_{B_r}|^2,
\end{aligned}\]
and hence
\[
 \mathcal E(u,B_r) 
 \le  2c_1\epsilon_0|(V(Du))_{B_r}|^2
 \le  \epsilon_*(\tau)|(V(Du))_{B_r}|^2.
\]
Thus Lemma~\ref{lem:small-excess} and the choice of $\tau$ yield
\[
 \mathcal E(u,B_{\tau r}) \le\frac14\mathcal E(u,B_r) .
\]
If $\mathcal E(u,B_r) > \epsilon_0M(r/2)$ but $M(r/4)> \delta M(r/2)$,
then \eqref{eq:sup-drop} and the choice of $\delta$ imply
\[\begin{aligned}
 \mathcal E(u,B_{\tau r}) 
 \le(4\tau)^{-n}\mathcal E(u,B_{r/4}) 
 &\le  c_1(4\tau)^{-n} \left(M(r/2)-M(r/4)\right)\\
 &\le c_1(4\tau)^{-n}(1-\delta)M(r/2)
 <\frac18\mathcal E(u,B_{ r}).
\end{aligned}\]
In the remaining case, we have
\[
\mathcal E(u,B_r)
 >\epsilon_0M(r/2)
 \quad \text{and}\quad 
M(r/4) \le\delta M(r/2).
\]
Since $\tau<1/4$, it follows that $M(\tau r) \le\delta M(r/2)$.
Moreover, using $t\varphi'(t)\le q\varphi(t)$, we obtain
\[
 \mathcal E(u,B_{\tau r})
 \le\fint_{B_{\tau r}}|V(Du)|^2\,dx
 \le qM(\tau r)
 \le q\delta M(r/2).
\]

For $j\ge0$, set
\[
 r_j:=\tau^j r,\qquad
 E_j:=\mathcal E(u,B_{r_j}),\qquad
 M_j:=M(r_j/2).
\]
Since $\tau<1/4$, the preceding alternatives give either
\[
 E_{j+1}\le\frac14E_j,
\]
or
\[
 E_j>\epsilon_0M_j,\qquad
 M_{j+1}\le\delta M_j,\qquad
 E_{j+1}\le q\delta M_j.
\]
Moreover, in both cases,
\[
 M_{j+1}\le M_j
\quad \text{and} \quad 
E_{j+1}\le qM_j.
\]
The standard iteration argument in the proof of
\cite[Theorem~3.1]{GiaquintaModica1986}
(see also \cite[Theorem~6.4]{DSV2009})
therefore yields constants
$c>0$ and $\kappa\in(0,1)$, depending only on
$n,N,p,q$, such that
\[
 \mathcal E(u,B_{\tau^j r})
 \le c\kappa^j\mathcal E(u,B_r)
 \qquad(j\ge0).
\] 
For $0<\rho<r$, choose $j\ge0$ such that
$\tau^{j+1}<\rho/r\le\tau^j$.
The minimizing property of the mean gives
\[
 \mathcal E(u,B_{\rho})
 \le\left(\frac{\tau^jr}{\rho}\right)^n
       \mathcal E(u,B_{\tau^j r})
 \le c\tau^{-n}\kappa^j\mathcal E(u,B_r).
\]
Choosing
\[
 \sigma:= \min\left\{\frac{\log\kappa}{\log\tau},\frac12\right\} \in(0,1),
\]
and using $\kappa^j\le\tau^{j\sigma}
\le\tau^{-\sigma}(\rho/r)^\sigma$, we conclude that
\[
 \mathcal E(u,B_\rho)
 \le c\left(\frac{\rho}{r}\right)^\sigma\mathcal E(u,B_r),
 \qquad 0<\rho<r.
\]
This proves \eqref{eq:main-excess}.

Finally, arguing as in the deduction of
\cite[Theorem~6.5]{DSV2009} from
\cite[Theorem~6.4 and Lemma~2.10]{DSV2009},
we obtain \eqref{eq:strong-main-holder} from
\eqref{eq:main-excess} and \eqref{eq:Lipschitz}.
Here the inverse estimate for $V$ follows from
Lemma~\ref{lem:secant} and \eqref{eq:flux-indices}, and we may take
\[
 \gamma=\frac{\sigma}{\max\{2,q\}}>0.
\]
Thus $u\in C^{1,\gamma}_{\mathrm{loc}}(\Omega;\R^N)$,
with all constants and exponents depending only on
$n,N,p,q$.
\end{proof}

\section{Higher gradient integrability and H\"older regularity under weak growth}\label{sec:weak-proof}

\subsection{Smooth approximation and comparison minimizers}

We first construct smooth, strictly convex approximations of $\varphi$ and establish the growth estimates needed for the comparison arguments. This construction is motivated by the approximation in \cite[Section~5]{HastoOk2022}.
Fix $\chi\in C_c^\infty((0,1))$ with $\chi\ge 0$, $\int_0^1\chi(s) \,ds=1$, and $\|\chi\|_\infty+\|\chi'\|_\infty\le c(n)$. For $\epsilon\in(0,1)$, define
\begin{equation}\label{phi_epsilon}
 \phi_\epsilon(t):=
 \psi_\epsilon(t)+\epsilon\int_0^t\frac{\psi_\epsilon(s)}s\,ds,
\quad \text{where }\ \psi_\epsilon(t):=
 \int_0^1\chi(s)\varphi((1+\epsilon s)t)\,ds.
\end{equation}
We immediately see that  $\phi(t)\le \psi_\epsilon(t)\le \phi_\epsilon(t)$ for $t\ge 0$.

\begin{lemma}\label{lem:approximation}
Assume that \eqref{eq:potential-scaling} holds. Then $\phi_\epsilon\in C^1([0,\infty))\cap C^\infty((0,\infty))$ and, for every $t>0$,
\begin{equation}\label{eq:relative-approximation}
 \varphi(t)\le\phi_\epsilon(t)\le
 (1+\epsilon)^q\left(1+\frac{\epsilon}{p}\right)\varphi(t) 
\end{equation}
\begin{equation}\label{eq:relative-approximation1}
 p\phi_\epsilon(t)\le t\phi_\epsilon'(t)\le q\phi_\epsilon(t),
\end{equation}
\begin{equation}\label{eq:approx-derivative-comparison}
 \epsilon(p-1)\varphi(t)
 \le t^2 \phi_\epsilon''(t)
 \le \frac{c}{\epsilon} \varphi(t),
\end{equation}
where $c\ge1$ is a constant depending only on $n,p,q$.
\end{lemma}

\begin{proof}
For $t>0$, the representation
\[
 \psi_\epsilon(t)=\frac1{\epsilon t}
 \int_0^\infty\chi\left(\frac{s-t}{\epsilon t}\right)\varphi(s)\,ds
\]
shows that it is smooth and so is $\phi_\epsilon$.
By \eqref{eq:potential-scaling} and \eqref{phi_epsilon}, $\psi_\epsilon(s)\le(1+\epsilon)^q\varphi(s)\le (1+\epsilon)^q (s/t)^p\varphi(t)$ for all $0<s\le t$, which implies \eqref{eq:relative-approximation}. 
Since 
\begin{equation}\label{lem:approximation_pf0}
 \psi_\epsilon'(t)=
 \int_0^1\chi(s)(1+\epsilon s)\varphi'((1+\epsilon s)t)\,ds,
\end{equation}
applying \eqref{eq:weak-indices} with $t$ replaced by $(1+\epsilon s)t$  gives
\begin{equation}\label{lem:approximation_pf1}
 p\psi_\epsilon(t)\le t\psi_\epsilon'(t)\le q\psi_\epsilon(t)
 \quad \text{and}\quad 
 \lambda^p\psi_\epsilon(t)\le \psi_\epsilon(\lambda t)\le \lambda^q \psi_\epsilon(t) \ \ \text{for all }\ \lambda\ge1,
\end{equation}
and hence 
\[
 \frac{\psi_\epsilon(t)}q
 = \int_0^t \frac{1}s \left(\frac st\right)^q\psi_\epsilon(t)\,ds
 \le\int_0^t\frac{\psi_\epsilon(s)}s\,ds
 \le\int_0^t\frac{1}s\left(\frac st\right)^p\psi_\epsilon(t)\,ds
 =\frac{\psi_\epsilon(t)}p.
\]
Since $\phi_\epsilon'(t)= \psi_\epsilon'(t)+\epsilon \psi_\epsilon(t)/t$,
these inequalities prove \eqref{eq:relative-approximation1}. 
Moreover, since 
\[
 0\le\phi_\epsilon'(t) \le(q+\epsilon)(1+\epsilon)^q\frac{\varphi(t)}t \longrightarrow0
 \quad\text{as }\ t\to 0^+,
\]
$\phi_\epsilon\in C^1([0,\infty))$ with $\phi_\epsilon'(0)=0$. 

Direct differentiation gives
\[
 \phi_\epsilon''(t)=\psi_\epsilon''(t)
 +\epsilon\frac{t\psi_\epsilon'(t)-\psi_\epsilon(t)}{t^2}.
\]
The convexity of $\varphi$ implies $\psi_\epsilon''(t)\ge0$ for $t>0$. Combining this with the first inequality in \eqref{lem:approximation_pf1}, we obtain
\[
 \phi_\epsilon''(t)
 \ge\epsilon(p-1)\frac{\psi_\epsilon(t)}{t^2}
 \ge\epsilon(p-1)\frac{\varphi(t)}{t^2}>0.
\]
This gives the lower bound of \eqref{eq:approx-derivative-comparison}. For the upper bound, integration by parts in the smoothing parameter gives
\[
 \psi_\epsilon''(t)
 =-\frac1{\epsilon t}\int_0^1
 \frac{d}{ds}\big[\chi(s)(1+\epsilon s)^2\big]
 \varphi'((1+\epsilon s)t)\,ds.
\]
This identity is valid by distributional differentiation even when $\varphi'$ is discontinuous. There are no boundary terms, since $\chi$ has compact support in $(0,1)$. Since
\[
 \varphi'((1+\epsilon s)t)
 \le q\frac{\varphi((1+\epsilon s)t)}{(1+\epsilon s)t}
 \le q(1+\epsilon)^q\frac{\varphi(t)}t \le  q2^q\frac{\varphi(t)}t  \le  q2^q\frac{\varphi_\epsilon(t)}t 
\]
by \eqref{eq:weak-indices} and $t\psi_\epsilon'(t)-\psi_\epsilon(t) \le (q-1)\psi_\epsilon(t)\le (q-1)\phi_\epsilon(t)$ by \eqref{lem:approximation_pf1}, we obtain
\[
 0\le\phi_\epsilon''(t)\le  \frac{q2^{q+2}}{\epsilon} (\|\chi\|_{\infty}+\|\chi'\|_\infty)\frac{\phi_\epsilon(t)}{t^2}+\epsilon (q-1)\frac{\phi_\epsilon(t)}{t^2} \le \frac{c}{\epsilon} \frac{\phi(t)}{t^2},
\]
where in the last inequality, we have used the second inequality in \eqref{eq:relative-approximation}.
\end{proof}

From the previous lemma, we see that $W^{1,\phi}(\Omega)=W^{1,\phi_\epsilon}(\Omega)$ and $W^{1,\phi}_0(\Omega)=W^{1,\phi_\epsilon}_0(\Omega)$ if $|\Omega|<\infty$.

\begin{lemma}\label{lem:common-energy-limit}
Assume that \eqref{eq:potential-scaling} holds. Let  $u\in W^{1,\varphi}_{\mathrm{loc}}(\Omega;\R^N)$ be a weak solution to \eqref{main_system}, and let $B_r\Subset\Omega$. Suppose that 
\[
 M:=\esssup_{B_r}|u|<\infty.
\]
Let  $\{\epsilon_j\}$ be a sequence such that $0<\epsilon_j<1$ and $\epsilon_j\to0$ as $j\to \infty$. For each $j\ge1$, let $w_j\in u+W^{1,\varphi}_0(B_r;\R^N)$ be the unique minimizer of
\[
 w\longmapsto\int_{B_r}\phi_{\epsilon_j}(|Dw|)\,dx,
\]
or equivalently the weak solution to $\mathrm{div}\left(\frac{\phi_{\epsilon_j}'(|Dw_j|)}{|Dw_j|}Dw_j\right)=0$ in $B_r$. Then 
\[
\esssup_{B_r}|w_j| \le M,
\]
and up to a subsequence,
$w_j\rightharpoonup v$ weakly in $W^{1,\varphi}(B_r;\R^N)$ for some $v\in u+W^{1,\varphi}_0(B_r;\R^N)$, with
\[
 \int_{B_r}\varphi(|Dv|)\,dx=\int_{B_r}\varphi(|Du|)\,dx.
\]
\end{lemma}

\begin{proof}
By a translation of the coordinates and a suitable choice of representatives of $u$ and its trace, we may assume that $|u(x)|\le M$ for all $x\in\overline{B_r(0)}$. 
Since $\phi_{\epsilon_j}$ satisfies \eqref{eq:relative-approximation1}, the convex hull property \cite[Lemma~4.1]{BDNS2025}, which is the vector-valued counterpart of the maximum principle, applies to $w_j$. Consequently, we obtain $\esssup_{B_r}|w_j|\le M$.

By \eqref{eq:relative-approximation} and the minimality of $w_j$,
\[\begin{aligned}
 \int_{B_r}\varphi(|Dw_j|)\,dx
 \le\int_{B_r}\phi_{\epsilon_j}(|Dw_j|)\,dx
& \le\int_{B_r}\phi_{\epsilon_j}(|Du|)\,dx\\
& \le(1+\epsilon_j)^q\left(1+\frac{\epsilon_j}{p}\right)\int_{B_r}\varphi(|Du|)\,dx.
\end{aligned}\]
Therefore,
$\{w_j\}$ is bounded in $W^{1,\varphi}(B_r;\R^N)$.
By reflexivity, there exists $v\in W^{1,\varphi}(B_r;\R^N)$
such that, up to a subsequence,
$w_j\rightharpoonup v$ in $W^{1,\varphi}(B_r;\R^N)$.
By weak lower semicontinuity, we have
\[
 \int_{B_r}\varphi(|Dv|)\,dx
 \le\liminf_{j\to\infty}\int_{B_r}\varphi(|Dw_j|)\,dx
 \le\int_{B_r}\varphi(|Du|)\,dx.
\]
For the reverse inequality, since
$u+W^{1,\varphi}_0(B_r;\R^N)$ is closed and convex in
$W^{1,\varphi}(B_r;\R^N)$, Mazur's theorem implies that
$v\in u+W^{1,\varphi}_0(B_r;\R^N)$.
The minimality of $u$ in $B_r$ then yields
\[
 \int_{B_r}\varphi(|Du|)\,dx
 \le\int_{B_r}\varphi(|Dv|)\,dx,
\]
which completes the proof.
\end{proof}

The following consequence of Theorem~\ref{thm:gradient}
will be used later to prove higher integrability of $Du$.

\begin{lemma}\label{cor:approximation-regularity}
Assume that  \eqref{eq:potential-scaling} holds, and let $\phi_\epsilon$
be given by Lemma~\ref{lem:approximation}.
For each fixed $\epsilon\in(0,1)$, every weak solution 
$w\in W^{1,\phi_\epsilon}(\Omega)$ of 
\[\diver\left(\frac{\phi'_\epsilon(|Dw|)}{|Dw|} Dw\right) =0
\qquad\text{in }\ \Omega
\]
satisfies $w\in C^{1,\gamma_\epsilon}_{\mathrm{loc}}(\Omega)$ for some $\gamma_\epsilon\in(0,1)$ and $w\in C^\infty\big(\{x\in \Omega: |Dw(x)|>0\}\big)$.
The exponent $\gamma_\epsilon$ and the regularity constants may depend on
$\epsilon$.
\end{lemma}

\begin{proof}
The derivative comparisons in
\eqref{eq:approx-derivative-comparison} imply
\[
 c_\epsilon^{-2}
 \le\frac{t\phi_\epsilon''(t)}{\phi_\epsilon'(t)}
 \le c_\epsilon^2
 \qquad \text{for }\ t>0.
\]
Thus $\phi_\epsilon$ satisfies \eqref{eq:flux-indices}
with indices $p=1+c_\epsilon^{-2}$ and $q=1+c_\epsilon^2$.
So Theorem~\ref{thm:gradient} gives
\[
 w\in C^{1,\gamma_\epsilon}_{\mathrm{loc}}(\Omega)
\]
for some $\gamma_\epsilon\in(0,1)$. No estimate uniform
in $\epsilon$ is needed here.

Set $U:=\{x\in\Omega:|Dw(x)|>0\}$. Since
$w\in C^{1,\gamma_\epsilon}_{\mathrm{loc}}(\Omega)$,
the set $U$ is open, and $|Dw|$ is bounded above and
away from zero on each compact subset of $U$.
The nonlinearity
$A_\epsilon(P):=\phi_\epsilon'(|P|)P/|P|$
is smooth for $P\ne0$, and the system is nondegenerate in $U$.
A standard difference-quotient argument gives
$w\in W^{2,2}_{\mathrm{loc}}(U)$, and each $\partial_k w$
satisfies
\[
 \diver\big(DA_\epsilon(Dw)D(\partial_k w)\big)=0
 \qquad\text{in }U
\]
in the weak sense.
This linear system has locally H\"older continuous
coefficients and is uniformly elliptic on each compact
subset of $U$. Interior Schauder estimates for elliptic
systems and a standard bootstrap argument therefore yield
$w\in C^\infty(U)$; see \cite[Chapter~III]{Giaquinta}.
\end{proof}

\subsection{Gradient estimates in $L^s$ and proof of Theorem~\ref{thm:main}}
We first establish a Caccioppoli inequality for the normalized
gradient of $C^1$ weak solutions and apply it to weak solutions of the
approximating systems, whose regularity follows from
Lemma~\ref{cor:approximation-regularity}.
We then derive higher integrability estimates that are
uniform in the approximation parameter and pass to the limit
to prove Theorem~\ref{thm:main} and Corollary~\ref{rem:weak-growth}.

Let $\Psi\in C^1([0,\infty))\cap C^\infty((0,\infty))$ be an N-function such that $\Psi'(t),\Psi''(t)>0$ for all $t>0$.
Let $w\in C^1(B_r;\R^N)$ be a weak solution of
\begin{equation}\label{eq:smooth-aux-system}
 \int_{B_r}\Psi'(|Dw|)\frac{Dw}{|Dw|}:D\zeta\,dx=0
 \quad\text{for all }\zeta\in C_c^\infty(B_r;\R^N).
\end{equation}
Suppose that $w$ is smooth on
\[
 U:=\{x\in B_r:|Dw(x)|>0\}.
\]
\begin{lemma}\label{lem:angular}
Under the above setting, let $\eta\in C_c^\infty(B_r)$ be nonnegative, and  $g\in C^1([0,\infty))$ be a nonnegative, nondecreasing function vanishing identically in a neighborhood of $0$. Then
\begin{equation}\label{eq:angular}
 \begin{aligned}
 \int_U\eta^2 g(|Dw|)|Dw|&\Psi'(|Dw|)\left|D\left[\frac{Dw}{|Dw|}\right]\right|^2\,dx\\
 &\le c\int_U g(|Dw|)|Dw|\Psi'(|Dw|)|D\eta|^2\,dx
\end{aligned}\end{equation}
for some $c>0$ depending only on $n$. 
Moreover, suppose that, for some $1<p\le q<\infty$,
\[
 t\Psi'(t)\le q\Psi(t)
 \qquad\text{and}\qquad
 t\Psi''(t)\ge(p-1)\Psi'(t)
 \qquad\text{for all }\ t>0.
\]
Then
\begin{equation}\label{eq:angular-full}
\begin{aligned}\int_U\eta^2\frac{\Psi'(|Dw|)}{|Dw|}
 \Big[g(|Dw|)|D^2w|^2 &+|Dw|g'(|Dw|)| D |Dw||^2\Big]dx \\
 &\le c\int_U\Psi(|Dw|)g(|Dw|)|D\eta|^2\,dx,
\end{aligned}\end{equation}
where $c>0$ depends only on $n,p$ and $q$.
\end{lemma}

\begin{proof}
On $U$, we set
\[
 \Theta=[\Theta^\alpha_i]:=\frac{Dw}{|Dw|}=\left[\frac{\partial_iw^\alpha}{|Dw|}\right], \qquad i=1,\dots,n \ \ \text{and}\ \ \alpha=1,\dots,N.
\]
Then, we have
\[
 |\Theta|=1,
 \quad 0\le\Theta^T\Theta\le I_n,
 \quad\text{and}\quad
 \sum_{i=1}^n \sum_{\alpha=1}^N\Theta^\alpha_i\partial_k\Theta^\alpha_i=0 \quad \text{for }\ k=1,\dots,n.
\]
Since $w$ is smooth in $U$, the system \eqref{eq:smooth-aux-system} takes the pointwise form
\begin{equation}\label{eq:pde-direction}
 \Psi''(|Dw|)\sum_{i}\Theta^\alpha_i\partial_i|Dw|
 =-\Psi'(|Dw|)\sum_{i}\partial_i\Theta^\alpha_i \quad 
 \quad\text{in }\ U,
 \end{equation}
 for each $\alpha =1,\dots,N$.
Moreover, differentiating the system \eqref{eq:smooth-aux-system} in $x_k$ gives, for each $\alpha=1,\dots,N$ and $k=1,\dots,n$,
\[
 \sum_i\partial_i\partial_k\big(\Psi'(|Dw|)\Theta^\alpha_i\big)=0
 \qquad\text{in }\ U.
\]
The function
$\partial_k w^\alpha g(|Dw|)\eta^2$ belongs to $C_c^1(U)$,
since $g$ vanishes near zero and $Dw$ is continuous.
By testing the previous equation with this function and 
summing over $\alpha,k$, we obtain
\begin{align}\label{eq:differentiated-identity}
 0&=\int_U\eta^2 g(|Dw|)\sum_{\alpha,i,k}\partial_k\big(\Psi'(|Dw|)\Theta^\alpha_i\big)\partial_i\partial_k w^\alpha\,dx\notag\\
 &\quad +\int_U\eta^2 g'(|Dw|)\sum_{\alpha,i,k}\partial_k\big(\Psi'(|Dw|)\Theta^\alpha_i\big)\partial_k w^\alpha\partial_i|Dw|\,dx\notag\\
 &\quad +2\int_U\eta g(|Dw|)\sum_{\alpha,i,k}\partial_k\big(\Psi'(|Dw|)\Theta^\alpha_i\big)\partial_k w^\alpha\partial_i\eta\,dx.
\end{align}
We first consider the first term. Since $\partial_i\partial_k w^\alpha=\partial_k\big(|Dw|\Theta_i^\alpha\big)= \partial_k|Dw| \Theta^\alpha_i + |Dw|  \partial_k \Theta^\alpha_i$,
\begin{align*}
 &\sum_{\alpha,i,k}\partial_k\big(\Psi'(|Dw|)\Theta^\alpha_i\big)
                    \partial_i\partial_k w^\alpha\\
 &= \sum_{\alpha,i,k} \big(\Psi'(|Dw|)\partial_k\Theta^\alpha_i+\Psi''(|Dw|)\partial_k|Dw| \Theta^\alpha_i\big)
                   ( \partial_k|Dw| \Theta^\alpha_i + |Dw|  \partial_k \Theta^\alpha_i)\\
&= \sum_{\alpha,i,k}  \Psi'(|Dw|)|Dw| (\partial_k\Theta^\alpha_i)^2 
+  \Psi''(|Dw|)\sum_{k}  (\partial_k|Dw|)^2 \sum_{\alpha,i} (\Theta^\alpha_i)^2 \\
&\qquad
+\sum_{k} \Psi'(|Dw|)\partial_k|Dw| \underbrace{\sum_{\alpha,i} \Theta^\alpha_i \partial_k\Theta^\alpha_i   }_{=0}\\
&\qquad +   \Psi''(|Dw|)|Dw|\sum_{k} \partial_k|Dw| \underbrace{\sum_{\alpha,i} \Theta^\alpha_i \partial_k\Theta^\alpha_i   }_{=0} \\
 & =|Dw|\Psi'(|Dw|)|D\Theta|^2+\Psi''(|Dw|)| D |Dw||^2,\\
 \end{align*}
 Since $
 (\partial_i|Dw|)\Theta_k^\alpha +|Dw|\partial_i\Theta_k^\alpha
 = \partial_i\partial_k w^\alpha =  \partial_k\partial_i w^\alpha
 =(\partial_k|Dw|)\Theta_i^\alpha +|Dw|\partial_k\Theta_i^\alpha,
$
multiplying by $\Theta_k^\alpha$ and summing over
$\alpha,k$, we obtain
\[
 \partial_i|Dw|
 +|Dw|\underbrace{\sum_{\alpha,k}
       \Theta_k^\alpha\partial_i\Theta_k^\alpha}_{=0}
 =
 \sum_{\alpha,k}\Theta_i^\alpha\Theta_k^\alpha
       \partial_k|Dw|
 +|Dw|\sum_{\alpha,k}\Theta_k^\alpha
       \partial_k\Theta_i^\alpha.
\]
Thus
\begin{equation}\label{eq:theta-derivative}
 |Dw|\sum_{\alpha,k}\Theta_k^\alpha
       \partial_k\Theta_i^\alpha
 =\partial_i|Dw|-(\Theta^T\Theta D |Dw|)_i.
\end{equation}

We now compute the second integrand in
\eqref{eq:differentiated-identity}. Using
$\partial_k w^\alpha=|Dw|\Theta_k^\alpha$ and
\eqref{eq:theta-derivative}, we find that
\begin{align*}
 &\sum_{\alpha,k} \partial_k\big(\Psi'(|Dw|)\Theta_i^\alpha\big)
 \partial_k w^\alpha\\
&\qquad =|Dw|\Psi''(|Dw|)\sum_{\alpha,k} \Theta_i^\alpha\Theta_k^\alpha\partial_k|Dw| +|Dw|\Psi'(|Dw|)\sum_{\alpha,k} \Theta_k^\alpha\partial_k\Theta_i^\alpha\\
&\qquad =|Dw|\Psi''(|Dw|)(\Theta^T\Theta D |Dw|)_i +\Psi'(|Dw|) \big((I_n-\Theta^T\Theta) D |Dw|\big)_i.
\end{align*}
Multiplying by $\partial_i|Dw|$ and summing over $i$ gives
\begin{align*}
 \sum_{\alpha,i,k}
 \partial_k\big(\Psi'(|Dw|)\Theta_i^\alpha\big)
 \partial_k w^\alpha\partial_i|Dw|
 &=|Dw|\Psi''(|Dw|)
       \big(\Theta^T\Theta D |Dw|\big)\cdot  D |Dw|\\
 &\hspace{-1cm}+\Psi'(|Dw|)
       \big[(I_n-\Theta^T\Theta) D |Dw|\big]
              \cdot  D |Dw|.
\end{align*}
Consequently, after factoring out $\eta^2$, the sum of
the first two integrands in
\eqref{eq:differentiated-identity} equals
\begin{align*}
 &g(|Dw|)|Dw|\Psi'(|Dw|)|D\Theta|^2\\
 &\quad+\big(g'(|Dw|)\Psi'(|Dw|)+ g(|Dw|)\Psi''(|Dw|)\big) \big[(I_n-\Theta^T\Theta) D |Dw|\big]\cdot  D |Dw|\\
 &\quad +\big(g(|Dw|)+|Dw|g'(|Dw|)\big)\Psi''(|Dw|)
       (\Theta^T\Theta D |Dw|)\cdot D |Dw|.
\end{align*}
Note that all three terms in the last expression are nonnegative,
since $g,g'\ge0$, $\Psi',\Psi''>0$, and
$0\le\Theta^T\Theta\le I_n$.

We now estimate the last integral in \eqref{eq:differentiated-identity}. 
It follows from $\partial_k w^\alpha = \Theta^\alpha_k |Dw|$ and \eqref{eq:pde-direction}  that
\begin{align*}
 \sum_{\alpha,k} \partial_k\big(\Psi'(|Dw|)\Theta_i^\alpha\big) \partial_k w^\alpha
 &=|Dw|\Psi'(|Dw|)\sum_{\alpha,k} \Theta_k^\alpha\partial_k\Theta_i^\alpha\\
 &\qquad+|Dw|\sum_\alpha\Theta_i^\alpha \left(\Psi''(|Dw|)\sum_k \Theta_k^\alpha\partial_k|Dw|\right)\\
 &=|Dw|\Psi'(|Dw|)\sum_{\alpha,k} \big(\Theta_k^\alpha\partial_k\Theta_i^\alpha -\Theta_i^\alpha\partial_k\Theta_k^\alpha\big),
\end{align*}
and hence
 \[
\bigg| \sum_{\alpha,k} \partial_k\big(\Psi'(|Dw|)\Theta_i^\alpha\big) \partial_k w^\alpha \bigg| \le c(n) |Dw|\Psi'(|Dw|) |D\Theta|
\]
for some $c(n)>0$.

Using the nonnegativity of the first two integrands
and then Young's inequality, we conclude that
\begin{align*}
 &\int_U\eta^2 g(|Dw|)|Dw|\Psi'(|Dw|)|D\Theta|^2\,dx \\
 &\qquad \le2c(n)\int_U\eta g(|Dw|)|Dw|\Psi'(|Dw|)|D\Theta||D\eta|\,dx\\
 &\qquad\le\frac12\int_U\eta^2 g(|Dw|)|Dw|\Psi'(|Dw|)|D\Theta|^2\,dx\\
 &\qquad\qquad+2c(n)^2\int_U g(|Dw|)|Dw|\Psi'(|Dw|)|D\eta|^2\,dx.
\end{align*}
Absorbing the first term on the right-hand side proves
\eqref{eq:angular}.

To prove \eqref{eq:angular-full}, we retain all three
nonnegative terms in the sum of the first two integrands
in \eqref{eq:differentiated-identity}.
Using the same estimate for the cutoff term and absorbing
one half of the term involving $|D\Theta|^2$, we obtain
\begin{align*}
 &\frac12\int_U\eta^2 g(|Dw|)|Dw|\Psi'(|Dw|)|D\Theta|^2\,dx\\
 &\quad+\int_U\eta^2\big(g'(|Dw|)\Psi'(|Dw|)+ g(|Dw|)\Psi''(|Dw|)\big)\\ 
 &\qquad\qquad\qquad\qquad\qquad \times \big[(I_n-\Theta^T\Theta) D |Dw|\big]\cdot D |Dw| \,dx\\
 &\quad+\int_U\eta^2 \big(g(|Dw|)+|Dw|g'(|Dw|)\big)\Psi''(|Dw|)\big[\Theta^T\Theta D |Dw|\big]\cdot D |Dw|\,dx\\
 &\qquad\le c(n)\int_U g(|Dw|)|Dw|\Psi'(|Dw|)|D\eta|^2\,dx.
\end{align*}
Set $\nu:=\min\{1,p-1\}>0$.
Since $0\le\Theta^T\Theta\le I_n$ and
$|Dw|\Psi''(|Dw|)\ge(p-1)\Psi'(|Dw|)$,
the left-hand side is bounded below by
\begin{align*}
 &\frac{\nu}{2}\int_U\eta^2
       \frac{\Psi'(|Dw|)}{|Dw|}g(|Dw|)
       \left[|Dw|^2|D\Theta|^2+| D |Dw||^2\right]dx\\
 &\quad+\nu\int_U\eta^2
       \Psi'(|Dw|)g'(|Dw|)| D |Dw||^2\,dx\\
 &\qquad\ge\frac{\nu}{2}\int_U\eta^2
       \frac{\Psi'(|Dw|)}{|Dw|}
       \left[g(|Dw|)|D^2w|^2
             +|Dw|g'(|Dw|)| D |Dw||^2\right]dx,
\end{align*}
where we used
\[
 |D^2w|^2
 =|Dw|^2|D\Theta|^2+| D |Dw||^2.
\]
Combining these inequalities and using
$|Dw|\Psi'(|Dw|)\le q\Psi(|Dw|)$ proves
\eqref{eq:angular-full}.
\end{proof}

We now use the Caccioppoli inequality for the normalized
gradient to obtain higher integrability estimates for
the gradient of bounded weak solutions.

\begin{proposition}\label{prop:moment}
Suppose that
$\Psi\in C^1([0,\infty))\cap C^\infty((0,\infty))$
is an N-function with $\Psi''(t)>0$ for
$t>0$ and satisfies
\[
 1<p\le\frac{t\Psi'(t)}{\Psi(t)}\le q <\infty
 \qquad\text{for all }\ t>0.
\]
Let $w\in C^1(B_r;\R^N)\cap C^\infty(U;\R^N)$ be a weak
solution to the system \eqref{eq:smooth-aux-system}
with
\[
 |w|\le M\quad\text{in }B_r
\]
for some $M>0$.
Then for every $b>\max\{1,(p-1)^{-1}\}$, there exists
a constant $c_b>0$, depending only on $n,p,q$ and $b$,
such that
\begin{equation}\label{eq:flux-moment}
 \fint_{B_{r/2}}|Dw|\Psi'(|Dw|)^b\,dx
 \le c_b\frac{M}{r}\Psi'\left(\frac{M}{r}\right)^b.
\end{equation}
Consequently, for every $s\in(1,\infty)$, there exists
$c_s>0$, depending only on $n,p,q$ and $s$, such that
\begin{equation}\label{eq:uniform-moment}
 \left(\fint_{B_{r/2}}|Dw|^s\,dx\right)^{1/s}
 \le c_s\frac{M}{r}.
\end{equation}
\end{proposition}

\begin{proof}
Fix $b>\max\{1,(p-1)^{-1}\}$ and, for $\kappa>0$, set
\[
 h_\kappa(t):=
 \big(t^{(b-1)/2}-\Psi'(\kappa)^{(b-1)/2}\big)_+^2,
 \quad t\ge0.
\]
Then, since $b>1$ and $\Psi'(\kappa)>0$, $h_\kappa$ is in $C^1([0,\infty))$ and increasing,
and both its value and derivative vanish at
$t=\Psi'(\kappa)$.
For $r/2\le s<\rho\le3r/4$, choose
$\eta\in C_c^\infty(B_\rho)$ such that $0\le\eta\le1$,
$\eta\equiv1$ on $B_s$, and $|D\eta|\le2/(\rho-s)$.
Hence the function
$w h_\kappa(\Psi'(|Dw|))\eta^2$
has compact support in $U$.
Testing \eqref{eq:smooth-aux-system} with this function
and using \eqref{eq:pde-direction} gives
\begin{equation}\label{eq:original-test}
\begin{aligned}
 &\int_U\eta^2|Dw|\Psi'(|Dw|)
       h_\kappa(\Psi'(|Dw|))\,dx\\
 &\qquad=-\int_U
       \eta^2\Psi'(|Dw|)h_\kappa'(\Psi'(|Dw|))
       \Psi''(|Dw|)
       \big[\Theta:w\otimes D(|Dw|)\big]\,dx\\
 &\qquad\qquad
       -2\int_U\eta\Psi'(|Dw|)
       h_\kappa(\Psi'(|Dw|))
       \big[\Theta:w\otimes D\eta\big]\,dx\\
 &\qquad=\int_U
       \eta^2\Psi'(|Dw|)^2h_\kappa'(\Psi'(|Dw|))
       \big[\diver\Theta\cdot w\big]\,dx\\
 &\quad\qquad
       -2\int_U\eta\Psi'(|Dw|)
       h_\kappa(\Psi'(|Dw|))
       \big[w\cdot\Theta D\eta\big]\,dx\\
 &\qquad=:I_1+I_2,
\end{aligned}
\end{equation}
where $\Theta:=Dw/|Dw|$.
Moreover, since
$h_\kappa'(t)^2
=h_\kappa(t)(b-1)^2t^{b-3}$ for $t>0$
and $|\diver\Theta|\le\sqrt n|D\Theta|$,
applying
\eqref{eq:angular} with
$g(t)=h_\kappa(\Psi'(t))\le\Psi'(t)^{b-1}$ and H\"older's inequality,
the absolute value of $I_1$ is estimated by
\begin{align*}
 |I_1|
 &\le\sqrt n M\int_U\eta^2\Psi'(|Dw|)^2 h_\kappa'(\Psi'(|Dw|))|D\Theta|\,dx\\
 &=\sqrt n M \int_{U\cap\{|Dw|>\kappa\}} \eta \sqrt{|Dw|\Psi'(|Dw|) h_\kappa(\Psi'(|Dw|))} \,|D\Theta|\\ 
 &\qquad\qquad\qquad\qquad\qquad\qquad\qquad \times \underbrace{\eta \sqrt{\frac{\Psi'(|Dw|)^3 h_\kappa'(\Psi'(|Dw|))^2} {|Dw|h_\kappa(\Psi'(|Dw|))}}}_{=(b-1)\eta\sqrt{\frac{\Psi'(|Dw|)^b}{|Dw|}}}\,dx\\
 &\le\sqrt n\,(b-1)M \left(\int_U \eta^2|Dw|\Psi'(|Dw|) h_\kappa(\Psi'(|Dw|))|D\Theta|^2\,dx \right)^{1/2}\\
 &\qquad\times \left(\int_{B_\rho}\eta^2\frac{\Psi'(|Dw|)^b}{|Dw|}\,dx \right)^{1/2}\\
 &\le\frac{c(n)(b-1)M}{\rho-s} \left(\int_{B_\rho} |Dw|\Psi'(|Dw|)^b\,dx\right)^{1/2} \left(\int_{B_\rho} \frac{\Psi'(|Dw|)^b}{|Dw|}\,dx\right)^{1/2}.
\end{align*}
For $I_2$, we use
$|\Theta|=1$, $|w|\le M$,
$h_\kappa(\Psi'(|Dw|))\le\Psi'(|Dw|)^{b-1}$, and
$|D\eta|\le2/(\rho-s)$

to obtain
\begin{align*}
 |I_2|
 &\le\frac{4M}{\rho-s}
       \int_{B_\rho}\Psi'(|Dw|)^b\,dx\\
 &\le\frac{4M}{\rho-s}
       \left(\int_{B_\rho}
       |Dw|\Psi'(|Dw|)^b\,dx\right)^{1/2}
       \left(\int_{B_\rho}
       \frac{\Psi'(|Dw|)^b}{|Dw|}\,dx\right)^{1/2}.
\end{align*}
Inserting the estimates for $I_1$ and $I_2$ into
\eqref{eq:original-test} and letting $\kappa\to0^+$,
by the monotone convergence theorem,
we obtain
\begin{equation}\label{eq:moment-key}
\begin{aligned} 
&\int_{B_s}|Dw|\Psi'(|Dw|)^b\,dx\\
&\qquad \le\frac{cbM}{\rho-s} \left(\int_{B_\rho} |Dw|\Psi'(|Dw|)^b\,dx\right)^{1/2} \left(\int_{B_\rho} \frac{\Psi'(|Dw|)^b}{|Dw|}\,dx\right)^{1/2},
\end{aligned}\end{equation}
where $c\ge1$ depends only on $n$.

We split $B_\rho$ into the sets where
$|Dw|<2cbM/(\rho-s)$ and $|Dw|\ge2cbM/(\rho-s)$.
By \eqref{eq:flux-scaling}, applied to $\Psi$,
and $b(p-1)>1$,
\[
 \frac{\Psi'(t)^b}{t}
 \le\left(\frac qp\right)^b\frac{\Psi'(s)^b}{s}
       \left(\frac ts\right)^{b(p-1)-1}
 \le\left(\frac qp\right)^b\frac{\Psi'(s)^b}{s}
 \qquad\text{for every }\ 0<t\le s.
\]
Hence, on $\{0<|Dw|<2cbM/(\rho-s)\}$,
using \eqref{eq:flux-scaling} once more,
\[
 \frac{\Psi'(\frac12|Dw|)^b}{\frac12|Dw|}
 \le\left(\frac qp\right)^b
       \frac{\Psi'(\frac{cbM}{\rho-s})^b}
            {\frac{cbM}{\rho-s}},
\]
and hence
\[
 \frac{\Psi'(|Dw|)^b}{|Dw|}
 \le c_b
       \left(\frac{r}{\rho-s}\right)^{(q-1)b}
       \frac{\rho-s}{M}
       \Psi'\left(\frac Mr\right)^b,
\]
where $c_b>0$ depends on $p,q,n$ and $b$.
Therefore, using this estimate in \eqref{eq:moment-key} and applying Young's inequality, we obtain
\[
\begin{aligned}
 &\int_{B_s}|Dw|\Psi'(|Dw|)^b\,dx\\
 &\le\frac12\int_{B_\rho}|Dw|\Psi'(|Dw|)^b\,dx\\
 &\qquad+c_b
 \left(\int_{B_\rho}|Dw|\Psi'(|Dw|)^b\,dx\right)^{1/2}
 \left( \left(\frac{r}{\rho-s}\right)^{(q-1)b+1}r^n\frac Mr\Psi'\left(\frac Mr\right)^b \right)^{1/2}\\
 &\le\frac34\int_{B_\rho}|Dw|\Psi'(|Dw|)^b\,dx+c_b\left(\frac{r}{\rho-s}\right)^{(q-1)b+1} r^n\frac Mr\Psi'\left(\frac Mr\right)^b.
\end{aligned}
\]
Since the integral over $B_{3r/4}$ is finite,
the standard iteration lemma
\cite[Lemma~3.1, p.~161]{Giaquinta} gives
\[
 \int_{B_{r/2}}|Dw|\Psi'(|Dw|)^b\,dx
 \le c_b r^n\frac Mr\Psi'\left(\frac Mr\right)^b.
\]
Dividing by $|B_{r/2}|$ proves \eqref{eq:flux-moment}.

The estimate \eqref{eq:uniform-moment} follows directly
from \eqref{eq:flux-moment} and \eqref{eq:flux-scaling}.
It is enough to consider $s>1+\max\{1,p-1\}$, since the estimates for smaller exponents follow by H\"older's inequality.
Then, choose $b$ such that
$s=1+(p-1)b$ and set
\[
 \Psi_{s,b}(t):=
 t^{1/s} \Psi'(t^{1/s})^b,
 \qquad t\ge0.
\]
Since
$\Psi_{s,b}(t)/t
 =\left(\Psi'(t^{1/s})/(t^{1/s})^{p-1}
 \right)^b$,
\eqref{eq:flux-scaling} shows that $\Psi_{s,b}(t)/t$ is almost increasing with constant $(q/p)^b$.
Moreover, for every $\lambda\ge1$ and $t>0$,
$\Psi_{s,b}(\lambda t) \le (q/p)^b\lambda^{(1+b(q-1))/s}\Psi_{s,b}(t)$.
By the Jensen-type inequality \cite[Lemma~4.3.1]{HarjulehtoHasto2019}, we obtain
\[\begin{aligned}
 \Psi_{s,b}\left(\fint_{B_{r/2}}|Dw|^s\,dx\right) 
&\le c_b\fint_{B_{r/2}} |Dw|\Psi'(|Dw|)^b\,dx \\
&\le c_b\frac Mr\Psi'\left(\frac Mr\right)^b=c_b \Psi_{s,b} \left(\left(\frac Mr\right)^s\right),
\end{aligned}\]
which, by the almost increasing property of
$\Psi_{s,b}(t)/t$, proves \eqref{eq:uniform-moment}.
\end{proof}

The limit of the approximating solutions has the same boundary values and energy as the original solution, but need not coincide with it.
The following lemma allows us to transfer the higher integrability estimates from this limit to the original solution.

\begin{lemma}\label{lem:flat-interval-transfer}
Assume that \eqref{eq:potential-scaling} holds, and let
$u\in W^{1,\varphi}(B_r;\R^N)$ be a weak solution
to \eqref{main_system} in $B_r$.
Suppose that $v\in u+W^{1,\varphi}_0(B_r;\R^N)$ satisfies
\[
 \int_{B_r}\varphi(|Dv|)\,dx
 =\int_{B_r}\varphi(|Du|)\,dx.
\]
Then
\begin{equation}\label{eq:transfer-gradient-comparison}
 |Du|\le\frac{p(q-1)}{q(p-1)}|Dv|
 \quad\text{almost everywhere in }\ B_r.
\end{equation}
\end{lemma}

\begin{proof}
Since $u$ minimizes
the $\varphi$-energy among functions with the same boundary value and $(u+v)/2$ has the same one,
the convexity of $\phi$ and the equality of the energies give
\begin{align*}
 \int_{B_r}\varphi(|Du|)\,dx
 &\le\int_{B_r} \varphi\left(\left|\frac{Du+Dv}{2}\right|\right)dx\\
 &\le\frac12\int_{B_r}  \big(\varphi(|Du|)+\varphi(|Dv|)\big)\,dx
 =\int_{B_r}\varphi(|Du|)\,dx.
\end{align*}
Moreover, pointwise convexity gives
\[
 \varphi\left(\left|\frac{Du+Dv}{2}\right|\right)
 \le\varphi\left(\frac{|Du|+|Dv|}{2}\right)
 \le\frac{\varphi(|Du|)+\varphi(|Dv|)}2.
\]
From the above results, we obtain
\begin{equation}\label{eq:lemv1v21}
 \varphi\left(\left|\frac{Du+Dv}{2}\right|\right)
 =\varphi\left(\frac{|Du|+|Dv|}{2}\right)
 =\frac{\varphi(|Du|)+\varphi(|Dv|)}2
 \quad\text{a.e. in }B_r.
\end{equation}

If $|Du|\le|Dv|$, then the conclusion follows from
$\frac{p(q-1)}{q(p-1)}\ge1$.
Suppose $|Du|>|Dv|$.
By the convexity of $\varphi$ and the second equality
in \eqref{eq:lemv1v21}, $\varphi$ is affine on
$[|Dv|,|Du|]$. Thus
$\varphi(t)=a_1t-a_0$ on this interval for some
$a_1>0$ and $a_0\in\R$.
Applying \eqref{eq:weak-indices} at $t\in(|Dv|,|Du|)$ gives
\[
 p\le\frac{a_1t}{a_1t-a_0}\le q,
 \quad\text{and hence}\quad
 \frac{p-1}{p}a_1t\le a_0
 \le\frac{q-1}{q}a_1t.
\]
Letting $t\nearrow|Du|$ in the first inequality and
$t\searrow|Dv|$ in the second, we obtain
\[
 \frac{p-1}{p}a_1|Du|\le a_0
 \le\frac{q-1}{q}a_1|Dv|,
\]
which proves \eqref{eq:transfer-gradient-comparison}.
\end{proof}

We now apply the higher integrability estimates to the
approximating solutions and pass to the limit.
Lemma~\ref{lem:flat-interval-transfer} then transfers these
estimates to the original solution.

\begin{proof}[Proof of Theorem~\ref{thm:main}]
Let $B_{2r}(x_0)\Subset\Omega$. The local boundedness theorem \cite[Theorem~1.1]{BDNS2025} yields
\[
 M:=\esssup_{B_r(x_0)}|u|<\infty.
\]
If $M=0$, the conclusions on this ball are immediate. 
Suppose $M>0$. Let $w_j$ and $v$ be as in Lemma~\ref{lem:common-energy-limit}.
By Lemma~\ref{cor:approximation-regularity}, each $w_j$ is $C^1$ in $B_r(x_0)$ and smooth on $\{|Dw_j|>0\}$.
By \eqref{eq:relative-approximation1} and \eqref{eq:approx-derivative-comparison}, the approximating functions $\phi_j:=\phi_{\epsilon_j}$ satisfy
\[
 p\le\frac{t\phi_j'(t)}{\phi_j(t)}\le q
 \quad \text{and}\quad
\phi_j''(t)>0
 \qquad\text{for all }\ t>0,
\]
with the same indices $p$ and $q$.
Since $\esssup_{B_r(x_0)}|w_j|\le M$, Proposition~\ref{prop:moment}, applied with $\Psi=\phi_j$, gives
\[
 \left(\fint_{B_{r/2}(x_0)}|Dw_j|^s\,dx\right)^{1/s}
 \le c_s\frac Mr
 \qquad\text{for every }s\in(1,\infty),
\]
where $c_s$ depends only on $n,p,q,s$ and is independent
of $j$. For each fixed $s\in(1,\infty)$, the sequence $\{Dw_j\}$
is bounded in $L^s(B_{r/2}(x_0);\R^{N\times n})$.
The already established weak $W^{1,\varphi}$ convergence
identifies every weak $L^s$ subsequential limit with $Dv$.
Hence $Dw_j\rightharpoonup Dv$ in
$L^s(B_{r/2}(x_0);\R^{N\times n})$.
Moreover, Lemma~\ref{lem:common-energy-limit} gives
\[
 v-u\in W^{1,\varphi}_0(B_r(x_0);\R^N),
 \qquad
 \int_{B_r(x_0)}\varphi(|Dv|)\,dx
 =\int_{B_r(x_0)}\varphi(|Du|)\,dx.
\]
We may therefore apply Lemma~\ref{lem:flat-interval-transfer}.
Together with weak lower semicontinuity, this yields
\begin{equation}\label{eq:main-gradient}
\begin{aligned}
 \left(\fint_{B_{r/2}(x_0)}|Du|^s\,dx\right)^{1/s}
 &\le\frac{p(q-1)}{q(p-1)}
       \left(\fint_{B_{r/2}(x_0)}|Dv|^s\,dx\right)^{1/s}\\
 &\le\frac{p(q-1)}{q(p-1)}\liminf_{j\to\infty}
       \left(\fint_{B_{r/2}(x_0)}|Dw_j|^s\,dx\right)^{1/s}
 \le c \frac Mr.
\end{aligned}\end{equation}
%

We now apply the preceding inequality to the translated solution $u-(u)_{B_{2r}(x_0)}$ and set
\[
 \widetilde{M}:=\esssup_{B_r(x_0)}|u-(u)_{B_{2r}(x_0)}|.
\]
The local boundedness theorem \cite[Theorem~1.1]{BDNS2025}, applied to the same solution, gives
\[
 \varphi\left(\frac{\widetilde{M}}{r}\right)
 \le c\fint_{B_{2r}(x_0)} \varphi\left(\frac{|u-(u)_{B_{2r}(x_0)}|}{r}\right)dx.
\]
By \eqref{eq:potential-scaling} and the Orlicz Sobolev--Poincar\'e inequality \cite[Lemma~A.1]{BDNS2025},
\[\begin{aligned}
 \fint_{B_{2r}(x_0)}\varphi\left(\frac{|u-(u)_{B_{2r}(x_0)}|}{r}\right)dx
& \le2^q\fint_{B_{2r}(x_0)} \varphi\left(\frac{|u-(u)_{B_{2r}(x_0)}|}{2r}\right)dx\\
&\le c\fint_{B_{2r}(x_0)}\varphi(|Du|)\,dx.
\end{aligned}\]
Taking the inverse and absorbing the constant by \eqref{eq:potential-scaling}, we obtain
\[
\frac{\widetilde{M}}r \le c\,\varphi^{-1}\left(\fint_{B_{2r}(x_0)}\varphi(|Du|)\,dx\right).
\]
Therefore, combining this with \eqref{eq:main-gradient} with $M$ replaced by $\widetilde{M}$ proves \eqref{eq:reverse-holder}.
\end{proof}

\begin{proof}[Proof of Corollary~\ref{rem:weak-growth}]
Set
\[
 p_*:=1+\frac1{2\{(2L)^{1/(p-1)}-1\}}
 \qquad\text{and}\qquad
 q_*:=L2^q-1.
\]
For $\lambda=(2L)^{1/(p-1)}>1$, convexity and
\eqref{eq:weak-growth} imply
\[
 t\varphi'(t) 
 \ge \frac{\varphi(t)-\varphi(t/\lambda)}{1-\lambda^{-1}}
 \ge \frac{1-L\lambda^{-p}}{1-\lambda^{-1}}\varphi(t)
 =p_*\varphi(t)
 \]
 and
 \[
 t\varphi'(t)
 \le\varphi(2t)-\varphi(t)
 \le(L2^q-1)\varphi(t)
 =q_*\varphi(t).
\]
Thus $\varphi$ satisfies \eqref{eq:weak-indices} with $1<p_*\le q_*<\infty$. Applying Theorem~\ref{thm:main} with these indices proves the claim.
\end{proof}

\section{Lipschitz regularity under a lower growth condition}\label{sec:lower-index-lipschitz}

In this section, we prove local
Lipschitz regularity without assuming an upper derivative bound. 
We first establish an estimate for smooth functions and solutions.
Its constant does not depend on an upper bound for the second derivative. We will then apply this estimate to the approximations
from Lemma~\ref{lem:approximation}.

\begin{lemma}\label{lem:lower-index-smooth}
Let $\Psi\in C^1([0,\infty))\cap C^\infty((0,\infty))$ be
an N-function satisfying
\begin{equation}\label{eq:lower-index-smooth-growth}
 p\le \frac{t\Psi'(t)}{\Psi(t)}\le q
 \qquad\text{and}\qquad
 \frac{t\Psi''(t)}{\Psi'(t)}\ge p-1
 \qquad\text{for all }\ t>0,
\end{equation}
where $1<p\le q<\infty$. Let
$w\in W^{1,\Psi}(B_{2r};\R^N)\cap C^1(B_{2r};\R^N)$ be
a weak solution of
\[
 \diver\left(\Psi'(|Dw|)\frac{Dw}{|Dw|}\right)=0
 \quad\text{in }B_{2r}.
\]
Suppose also that
$w\in C^\infty(U;\R^N)$, where
$U:=\{x\in B_{2r}:|Dw(x)|>0\}$. Then
\begin{equation}\label{eq:lower-index-smooth-bound}
 \esssup_{B_r}\Psi(|Dw|)  \le c\fint_{B_{2r}}\Psi(|Dw|)\,dx,
\end{equation}
where $c\ge1$ depends only on $n,N,p,q$.
\end{lemma}

\begin{proof}
\textit{Step 1: Powers of the energy.}
Fix $\ell\ge1$ and $\kappa>0$, and choose
\begin{equation}\label{eq:lower-index-weight-bounds}
 g_{\kappa,\ell}(t):= \Psi(t)^{\ell-1} \left(1-\frac{\Psi(\kappa)}{\Psi(t)}\right)_+^2 \le\Psi(t)^{\ell-1} \ \ \text{for }\  t>0
 \quad\text{and}\quad
 g_{\kappa,\ell}(0):=0.
\end{equation}
This function belongs to $C^1([0,\infty))$, is nonnegative
and nondecreasing, and vanishes on $[0,\kappa]$. 
The lower growth bound for $\Psi$ in \eqref{eq:lower-index-smooth-growth} also gives
\begin{equation}\label{eq:lower-index-weight-lower}
 g_{\kappa,\ell}(t)
 \ge(1-2^{-p})^2\Psi(t)^{\ell-1}
 \qquad\text{for }\ t>2\kappa.
\end{equation}

Set
\[
 \theta_\kappa(x):=\left(
       \Psi(|Dw(x)|)^{\ell/2}-\Psi(2\kappa)^{\ell/2}
       \right)_+, \qquad x\in B_{2r}.
\]
For each fixed $\kappa$, this function is locally Lipschitz.
On $\{|Dw|>2\kappa\}$, we have
\begin{align*}
 |D\theta_\kappa|^2
 &=\frac{\ell^2}{4}\Psi(|Dw|)^{\ell-2} \Psi'(|Dw|)^2| D |Dw||^2\\
 &\le\frac{q\ell^2}{4}\frac{\Psi'(|Dw|)}{|Dw|}\Psi(|Dw|)^{\ell-1}| D |Dw||^2\\
 &\le c\ell^2\frac{\Psi'(|Dw|)}{|Dw|} g_{\kappa,\ell}(|Dw|)|D^2w|^2,
\end{align*}
where we used \eqref{eq:lower-index-weight-lower} and
$| D |Dw||\le|D^2w|$.
Moreover, $D\theta_\kappa=0$ almost everywhere on the remaining set.
By \eqref{eq:lower-index-smooth-growth}, we may apply
\eqref{eq:angular-full} from Lemma~\ref{lem:angular}
on $B_{2r}$, in place of $B_{r}$, with $g=g_{\kappa,\ell}$.
Together with \eqref{eq:lower-index-weight-bounds}, this gives
\begin{equation}\label{eq:lower-index-energy-power0}
 \int_{B_{2r}}\eta^2|D\theta_\kappa|^2\,dx
 \le c\ell^2\int_{B_{2r}}\Psi(|Dw|)^\ell|D\eta|^2\,dx
\end{equation}
for every $\eta\in C_c^\infty(B_{2r})$ with $0\le\eta\le1$.
Since $Dw$ is continuous, $\theta_\kappa\to\Psi(|Dw|)^{\ell/2}$ locally uniformly
as $\kappa\to0^+$.
For each $B_\rho\Subset B_{2r}$, choose
$\zeta\in C_c^\infty(B_{2r})$ with $0\le\zeta\le1$ and
$\zeta\equiv1$ on $B_\rho$. The preceding estimate gives
\[
 \int_{B_\rho}|D\theta_\kappa|^2\,dx
 \le c\ell^2\int_{B_{2r}}\Psi(|Dw|)^\ell|D\zeta|^2\,dx,
\]
uniformly in $\kappa$.
This implies
that $\{\theta_\kappa\}$ is bounded in $W^{1,2}(B_\rho)$.
Thus, up to a subsequence,
$\theta_\kappa\rightharpoonup\Psi(|Dw|)^{\ell/2}$ in
$W^{1,2}(B_\rho)$, and hence
$\Psi(|Dw|)^{\ell/2}\in W^{1,2}_{\mathrm{loc}}(B_{2r})$.
Using weak
lower semicontinuity and \eqref{eq:lower-index-energy-power0}, we obtain that for every $\eta\in C_c^\infty(B_{2r})$ with $0\le\eta\le1$,
\begin{equation}\label{eq:lower-index-energy-power}
 \begin{aligned}
 \int_{B_{2r}}\eta^2 \left|D\big(\Psi(|Dw|)^{\ell/2}\big)\right|^2\,dx
 &\le\liminf_{\kappa\to0^+} \int_{B_{2r}}\eta^2|D\theta_\kappa|^2\,dx\\
 &\le c\ell^2\int_{B_{2r}}\Psi(|Dw|)^\ell|D\eta|^2\,dx.
\end{aligned}\end{equation}

\medskip
\textit{Step 2: Iteration.}
By translation and the change of variables
$\widetilde w(y)=w(x_0+2ry)/(2r)$, it is enough to work on $B_1$.
Set
\[
 \chi_0:=\frac{n}{n-2}\quad\text{if }\ n>2
 \qquad\text{and}\qquad
 \chi_0:=2\quad\text{if }\ n=2.
\]
For $1/2\le s<\rho\le3/4$, choose $\eta\in C_c^\infty(B_{\rho})$ with $0\le\eta\le1$,
$\eta\equiv1$ on $B_s$, and
$|D\eta|\le c/(\rho-s)$.
The Sobolev inequality, applied to $\eta \Psi(|Dw|)^{\ell/2}$,
and \eqref{eq:lower-index-energy-power} yield
\[
 \|\Psi(|Dw|)\|_{L^{\ell\chi_0}(B_s)}
 \le\left(\frac{c\ell^2}{(\rho-s)^2}\right)^{1/\ell}
       \|\Psi(|Dw|)\|_{L^\ell(B_\rho)}.
\]
Take $\ell_j:=\chi_0^j$ and
$r_j:=1/2+2^{-j-2}$ for $j\ge0$. Applying the last estimate
with $\ell=\ell_j$, $\rho=r_j$ and $s=r_{j+1}$, we get
\[
 \|\Psi(|Dw|)\|_{L^{\ell_{j+1}}(B_{r_{j+1}})}
 \le\left(c\chi_0^{2j}2^{2j+6}\right)^{\chi_0^{-j}}
       \|\Psi(|Dw|)\|_{L^{\ell_j}(B_{r_j})}.
\]
Since $\sum_{j=0}^\infty(1+j)\chi_0^{-j}<\infty$,
iteration gives
\[
 \esssup_{B_{1/2}}\Psi(|Dw|)  \le c\|\Psi(|Dw|)\|_{L^1(B_{3/4})}
 \le c\fint_{B_1}\Psi(|Dw|)\,dx.
\]
Rescaling proves \eqref{eq:lower-index-smooth-bound}.
\end{proof}

We now verify that the approximation preserves the condition
in \eqref{eq:lower-derivative-index} and pass to the limit in
the preceding estimate.

\begin{proof}[Proof of Theorem~\ref{thm:lower-index-lipschitz}]
Let $\psi_\epsilon$ and $\phi_\epsilon$ be as in
\eqref{phi_epsilon}. By \eqref{lem:approximation_pf0} and
\eqref{eq:lower-derivative-index}, for every $\lambda\ge1$
and $t>0$,
\begin{align*}
 \psi_\epsilon'(\lambda t)
 &=\int_0^1\chi(s)(1+\epsilon s)
       \varphi'((1+\epsilon s)\lambda t)\,ds\\
 &\ge\lambda^{p-1}\int_0^1\chi(s)(1+\epsilon s)
                 \varphi'((1+\epsilon s)t)\,ds
  =\lambda^{p-1}\psi_\epsilon'(t).
\end{align*}
Moreover, \eqref{eq:potential-scaling} and the definition of $\psi_\epsilon$ imply $\psi_\epsilon(\lambda t) \ge\lambda^p\psi_\epsilon(t)$.
Since $\phi_\epsilon'(t)=\psi_\epsilon'(t) +\epsilon\psi_\epsilon(t)/t$, we obtain
\begin{align*}
\phi_\epsilon'(\lambda t)
=\psi_\epsilon'(\lambda t) +\epsilon\frac{\psi_\epsilon(\lambda t)}{\lambda t}
\ge\lambda^{p-1}\left(\psi_\epsilon'(t)+\epsilon\frac{\psi_\epsilon(t)}t\right)
=\lambda^{p-1}\phi_\epsilon'(t).
\end{align*}
Thus $\phi_\epsilon'$ also satisfies $(\mathrm{Inc})_{p-1}$.
Since $\phi_\epsilon$ is smooth on $(0,\infty)$, we also have
\[
 t\phi_\epsilon''(t)\ge(p-1)\phi_\epsilon'(t)
 \quad\text{for all }t>0.
\]
Thus \eqref{eq:relative-approximation1} and this inequality
show that $\phi_\epsilon$ satisfies
\eqref{eq:lower-index-smooth-growth} with the same $p$ and $q$.

Fix $B_{2r}(x_0)\Subset\Omega$. By
\cite[Theorem~1.1]{BDNS2025}, $u$ is bounded on this ball.
Choose $\epsilon_j\in(0,1)$ with $\epsilon_j\to0$ as $j\to \infty$ and let $w_j$ be the weak solutions
with boundary values $u$ from Lemma~\ref{lem:common-energy-limit},
applied on $B_{2r}(x_0)$ in place of $B_r$.
By Lemma~\ref{cor:approximation-regularity},
$w_j\in C^1(B_{2r}(x_0);\R^N)$ and $w_j$ is smooth on
$\{|Dw_j|>0\}$. The upper derivative bounds in
\eqref{eq:approx-derivative-comparison} are used here only for this qualitative regularity; their constants may depend on $j$.
By Lemma~\ref{lem:lower-index-smooth}, the minimality of $w_j$,
and \eqref{eq:relative-approximation}, we have 
\begin{equation}\label{eq:lower-index-uniform}\begin{aligned}
 \esssup_{B_r(x_0)}\varphi(|Dw_j|)
 &\le\esssup_{B_r(x_0)}\phi_{\epsilon_j}(|Dw_j|)
 \le c\fint_{B_{2r}(x_0)}\phi_{\epsilon_j}(|Dw_j|)\,dx\\
 &\le c\fint_{B_{2r}(x_0)}\phi_{\epsilon_j}(|Du|)\,dx
 \le c \fint_{B_{2r}(x_0)}\varphi(|Du|)\,dx,
 \end{aligned}\end{equation}
where $c$ is independent of $j$.

By Lemma~\ref{lem:common-energy-limit}, up to a subsequence,
$w_j\rightharpoonup v$ in $W^{1,\varphi}(B_{2r}(x_0);\R^N)$,
where $v-u\in W^{1,\varphi}_0(B_{2r}(x_0);\R^N)$ and
\[
 \int_{B_{2r}(x_0)}\varphi(|Dv|)\,dx
 =\int_{B_{2r}(x_0)}\varphi(|Du|)\,dx.
\]
Since $\varphi'$ is strictly
increasing, $P\mapsto\varphi(|P|)$ is strictly convex.
By the uniqueness of the minimizer of the energy $\int_{B_{2r}(x_0)}\phi(|D\tilde u|)\,dx$ for $\tilde u \in u+W^{1,\phi}_0(B_{2r}(x_0))$, we obtain that $u=v$ almost everywhere in $B_{2r}(x_0)$.

By \eqref{eq:lower-index-uniform}, $\{Dw_j\}$ is uniformly bounded in $L^\infty(B_r(x_0);\R^{N\times n})$.
The same bound holds for $Du$ by the weak convergence in $W^{1,\varphi}$.
Together with the continuity of $\varphi$, this yields \eqref{eq:lower-index-lipschitz} and completes the proof.
\end{proof}

\section*{Acknowledgments}
The author was supported by the National Research Foundation
of Korea (NRF) grants funded by the Korean government (MSIT)
(NRF-2022R1C1C1004523 and RS-2025-24533680).

\medskip
\noindent\textbf{Data availability.}
No datasets were generated or analyzed during this study.

\medskip
\noindent\textbf{Conflict of interest.}
The author declares no conflict of interest.

\medskip
\noindent\textbf{AI statement.}
The author used ChatGPT to assist in checking intermediate calculations, clarifying and refining parts of the proofs,
and improving the English and presentation of the manuscript. 
The author takes full responsibility for the mathematical content and the final version of the manuscript.

\end{document}